\documentclass[a4paper,12pt]{amsart}

\numberwithin{equation}{section}
\theoremstyle{plain}
\newtheorem{theorem}{\sc \bf Theorem}[section]
\newtheorem{lemma}[theorem]{\sc \bf Lemma}
\newtheorem{corollary}[theorem]{\sc \bf Corollary}
\newtheorem{proposition}[theorem]{\sc \bf Proposition}
\newtheorem{claim}{\sc \bf Claim}
\theoremstyle{definition} 
\newtheorem{definition}[theorem]{\sc \bf Definition}
\newtheorem{remark}[theorem]{\sc \bf Remark}
\newtheorem{example}[theorem]{\sc \bf Example}

\usepackage{amssymb}
\usepackage{amsfonts}
\usepackage{amsmath}
\usepackage{mathrsfs}   
\usepackage{graphicx}
\usepackage{latexsym}
\usepackage{accents}
\usepackage{enumerate}
\usepackage{ulem}
\usepackage{fancybox}
\usepackage{wasysym}
\usepackage{color}
\usepackage{cancel}

\usepackage{tikz}

\newcommand{\K}{{\mathbb{K}}}

\newcommand{\R}{{\mathbb{R}}}
\newcommand{\C}{{\mathbb{C}}}

\newcommand{\ali}[1]{\begin{align*}#1\end{align*}}
\newcommand{\alil}[1]{\begin{align}#1\end{align}}

\newcommand{\Si}{\Sigma}
\newcommand{\si}{\sigma}

\newcommand{\ph}{\varphi}
\newcommand{\phe}{\varphi(\epsilon,\cdot)}
\newcommand{\lam}{\lambda}
\newcommand{\e}{\epsilon}
\newcommand{\p}{{\prime}}

\newcommand{\BR}{\mathcal{B}}

\newcommand{\DR}{\mathcal{D}}
\newcommand{\ER}{\mathcal{E}}

\newcommand{\IR}{\mathcal{I}}

\newcommand{\LR}{\mathcal{L}}

\newcommand{\OR}{\mathcal{O}}
\newcommand{\PR}{\mathcal{P}}
\newcommand{\QR}{\mathcal{Q}}
\newcommand{\RR}{\mathcal{R}}
\newcommand{\SR}{\mathcal{S}}
\newcommand{\TR}{\mathcal{T}}

\newcommand{\tLR}{\tilde{\mathcal{L}}}

\newcommand{\ite}[1]{\begin{enumerate}[(1)]#1\end{enumerate}}
\let\prop=\relax

\newcommand{\prop}[1]{\begin{proposition}#1\end{proposition}}
\newcommand{\props}{\begin{proposition}}
\newcommand{\prope}{\end{proposition}}

\newcommand{\cors}{\begin{corollary}}
\newcommand{\core}{\end{corollary}}
\newcommand{\thm}[1]{\begin{theorem}#1\end{theorem}}
\newcommand{\thms}{\begin{theorem}}
\newcommand{\thme}{\end{theorem}}

\newcommand{\lems}{\begin{lemma}}
\newcommand{\leme}{\end{lemma}}

\newcommand{\defis}{\begin{definition}}
\newcommand{\defie}{\end{definition}}

\newcommand{\exams}{\begin{example}}
\newcommand{\exame}{\end{example}}
\newcommand{\rem}[1]{\begin{remark}\normalfont #1\end{remark}}
\newcommand{\cla}[1]{\begin{claim}#1\end{claim}}
\newcommand{\ncla}[1]{\setcounter{claim}{0}\begin{claim}#1\end{claim}}

\newcommand{\pros}{\begin{proof}}
\newcommand{\proe}{\end{proof}}
\newcommand{\prossq}{\begin{proof}}
\newcommand{\case}[1]{\begin{cases}#1\end{cases}}
\newcommand{\cd}{\cdot}
\newcommand{\up}{\upsilon}
\newcommand{\om}{\omega}
\newcommand{\ti}[1]{\tilde{#1}}

\newcommand{\diam}{{\mathrm{diam}}}
\newcommand{\supp}{\mathrm{supp}}

\newcommand{\qqqquad}{\qquad\qquad}

\newcounter{constants}
\makeatletter
\def\addconst{
\addtocounter{constants}{1}
\def\@currentlabel{\arabic{constants}}
\@currentlabel
}
\makeatother

\newcommand{\adl}[1]{\addconst\label{c:#1}}
\newcommand{\adr}[1]{\ref{c:#1}}

\newcommand{\var}{\mathrm{var}}

\newcommand{\MatII}[1]{\left[\begin{array}{cc}#1\end{array}\right]}

\begin{document}
\title[Quasi-compactness of transfer operators]{Quasi-compactness of transfer operators for topological Markov shifts with holes}
\address{
{\rm Haruyoshi Tanaka}\\
Department of Mathematics and Statistics\\
Wakayama Medical University, 580, Mikazura, Wakayama-city, Wakayama, 641-0011
}
\email{htanaka@wakayama-med.ac.jp}
\author{Haruyoshi Tanaka}
\begin{abstract}
We consider transfer operators for topological Markov shift (TMS) with countable states and with holes which are $2$-cylinders. As main results, if the closed system of the shift has irreducible transition matrix and the potential is a weaker Lipschitz continuous and summable, then we obtain a version of Ruelle-Perron-Frobenius Theorem and quasi-compactness of the associated Ruelle transfer operator. The escape rate of the open system is also calculated. 
In corollary, it turns out that the Ruelle operator of summable potential on topologically transitive TMS has a spectral gap property. 
As other example, we apply the main results to the transfer operators associated to graph iterated function systems.

\end{abstract}
\maketitle
\setcounter{tocdepth}{3}
\section{Introduction and outline of main results}\label{sec:intro}
Let $S$ be a countable set and $A=(A(ij))_{S\times S}$ a zero-one matrix.
Consider the topological Markov shift $X=X_{A}$ with state space $S$ and transition matrix $A$, namely $X=\{\om=\om_{0}\om_{1}\cdots\in\prod_{n=0}^{\infty}S\,:\,A(\om_{i}\om_{i+1})=1 \text{ for any }i\geq 0\}$ with the shift transformation $\si\,:\,X\to X$ which is defined by $(\si\om)_{n}=\om_{n+1}$ for any $n\geq 0$.
We call a function $\ph\,:\,X\to \R$ {\it summable} if
\alil
{
\sum_{s\in S\,:\,[s]\neq \emptyset}\exp(\sup_{\om\in [s]}\ph(\om))<\infty, \label{eq:sum}
}
where for word $w\in S^{n}$, $[w]$ means the cylinder set $[w]=\{\om\in X\,:\,\om_{0}\cdots\om_{n-1}=w\}$. Such a summability is treated by \cite{MSU, MU, PU, RSU, SU} in mainly fractal analysis.
Note the difference of the terminologies `summable' and `summable variation'.
Here a function $\ph\,:\,X\to \R$ is said to be summable variation if
$\sum_{n=2}^{\infty}\var_{n}\ph<\infty$, where $\var_{n}\ph=\sup\{|\ph(\om)-\ph(\up)|\,:\,\om_{i}=\up_{i}\ \ (i=0,1,\dots, n-1)\}$.
For $\theta\in (0,1)$, a metric $d_{\theta}$ on $X$ is given by $d_{\theta}(\om,\up)=\theta^{\min\{n\geq 0\,:\,\om_{n}\neq \up_{n}\}}$ if $\om\neq \up$ and $d_{\theta}(\om,\up)=0$ if $\om=\up$.
Let $\K=\R$ or $\C$. For $f\,:\,X\to \K$ and $k\geq 1$, we define
\ali
{
[f]_{k}=\sup\{\var_{n}f/\theta^{n}\,:\,n\geq k\}.
}
Notice $[f]_{k}\geq [f]_{k+1}$. When $[f]_{1}<\infty$, $f$ is called {\it locally Lipschitz continuous}, and if $[f]_{2}<\infty$ then $f$ is called {\it weak Lipschitz continuous}.
Let $C(X,\K)$ be the set of all $\K$-valued continuous functions on $X$, and $F^{k}(X,\K)$ the set of all $f\in C(X,\K)$ with $[f]_{k}<\infty$. We define $C_{b}(X,\K)$ by the Banach space consisting of all bounded functions $f\in C(X,\K)$ endowed with the supremum norm $\|f\|_{\infty}=\sup_{X}|f|$ and $F^{k}_{b}(X,\K)$ by the Banach space consisting of all bounded functions $f\in F^{k}(X,\K)$ endowed with the norm $\|f\|_{k}=\|f\|_{\infty}+[f]_{k}$. For simplicity, `$\K$' is omitted from these definitions when $\K=\C$.
\smallskip
\par
Let $S_{0}\subset S$ be a nonempty subset and $M=(M(ij))$ a zero-one matrix indexed by $S_{0}$ such that $M(ij)\leq A(ij)$ for any $i,j\in S_{0}$. Assume the following three conditions:
\ite
{
\item[(A.1)] The matrix $A$ is irreducible. 
\item[(A.2)] The subsystem $X_{M}$ of $X_{A}$ with the transition matrix $M$ has a periodic point for $\si$.
\item[(A.3)] A function $\ph\,:\,X\to \R$ is summable and satisfies $[\ph]_{k+1}<+\infty$ for some $k\geq 1$.
}
(see Section \ref{sec:pre} for terminology). 
We introduce a transfer operator associated with an open system of the shift. An operator $\LR_{M}=\LR_{M,\ph}$ associated to $M$ and $\ph$ is defined by
\alil
{
\LR_{M} f(\om)=\LR_{M,\ph} f(\om):=\sum_{a\in S\,:\,M(a\om_{0})=1}e^{\ph(a\cd \om)}f(a\cd\om)\label{eq:transfer}
}
for $f\,:\,X\to \C$ and $\om\in X$ formally, where $a\cd \om$ means the concatenation $a\om_{0}\om_{1}\cdots$. 
Note that the operator $\LR_{M}$ acting on $F^{k}_{b}(X)$ and on $C_{b}(X)$ are both bounded. 
Such an operator is used in studying system with hole \cite{Demers_,T2020, T2019, T2009}. In fact when put $\Si=\bigcup_{ij\,:\,M(ij)=1}[ij]$, we regard the map $\si|_{\Si}\,:\,\Si\to X$ as an {\it open system} and the map $\si\,:\,X\to X$ as a {\it closed system}. 
Denoted by $P(\ph)$ the topological pressure of $\ph$ (see (\ref{eq:toppres})).
\smallskip
\par
Outline of our main results are the following (I)-(IV) under the conditions (A.1)-(A.3): 
\ite
{
\item[(I)] letting $\lam=\exp(P(\ph|_{X_{M}}))$, there exist a nonnegative function $g\in F^{k}_{b}(X,\R)$ with $\|g\|_{\infty}=1$ and a Borel probability measure $\nu$ on $X$ such that $\LR_{M}g=\lam g$ and $\LR_{M}^{*}\nu=\lam\nu$. In particular, $\lam$ equals the spectral radius of $\LR_{M}$ (Theorem  \ref{th:ex_efunc_geneRuelleop_M}(1));
\item[(II)] the essential spectral radius $r_{\mathrm{ess}}(\LR_{M})$ of $\LR_{M}$ acting on $F^{k}_{b}(X)$ satisfies $\theta\lam\leq r_{\mathrm{ess}}(\LR_{M})<\lam$ if $X_{A}$ is not reduced to finite single orbit, and equals $0$ otherwise (Theorem  \ref{th:ex_efunc_geneRuelleop_M}(2));
\item[(III)] the escape rate of the open system $\si|_{\Si}$ is $P(\ph|_{X_{M}})-P(\ph)$ (Theorem \ref{th:escaperate}); and
\item[(IV)] if $M$ is irreducible, then the peripheral eigenvalues of $\LR_{M}$ are the $p$-th roots of $z^{p}=\lam^{p}$ and simple, where $p$ is the period of $M$ (Theorem \ref{th:ex_efunc_geneRuelleop_M_irre}).
}
\par
In particular, if $A$ is finitely irreducible, then the essential spectral radius $r_{\mathrm{ess}}(\LR_{M})$ coincides with $\theta\lam$ exactly (Corollary \ref{cor:ex_efunc_geneRuelleop_M_fi}). These results are generalizations of \cite[Theorem 2.1]{Demers_} who gave the spectral gap property and the escape rate of $\LR_{M}$ under finite primitivity of $M$ and locally Lipschitz potential with finite pressure (therefore summable from \cite[Proposition 2.1.9]{MU}). 
When $A=M$, our results are also extensions of previous work in \cite{AD,MU2001} who considered under the finite irreducibility of $M$ and locally Lipschitz continuous of $\ph\,:\,X_{M}\to \R$. Moreover, \cite{Cyr_Sarig} introduced the notion of strongly positive recurrent which is a necessary and sufficient condition for the existence of a Banach space $F$ such that $\LR\,:\,F\to F$ is a bounded linear operator with spectral gap. Since $F$ may not contain a subset of $F^{k}_{b}(X_{M})$ (see Remark \ref{rem:SPR}), our result is different from \cite{Cyr_Sarig}. 
In \cite{T2009}, we gave similar spectral results for $\LR_{M}$ under the finite state case.
\smallskip
\indent
To prove quasi-compactness of Ruelle operators, the property of Gibbs measures was
used in previous study \cite{AD,MU}. Since this property imposes a topological restricting to the transition matrix concerning finite irreducibility, we need another technique in general TMS. The main tool of the proof of our results is a perturbation method of transfer operators involving a change of symbolic dynamics. Precisely, we introduce a perturbed system $(X_{A},\phe)$ such that $A$ is finitely irreducible, $\phe\in F^{1}(X_{A},\R)$ is summable uniformly in $\e>0$ and the operator $\LR_{A,\phe}$ converges to $\LR_{M}$ in a supremum norm 
(see (\ref{eq:ph-(1/e)chiN}) in Section \ref{sec:proof} and the conditions (B.1)-(B.3) in Section \ref{sec:aux}). By calculating the essential spectral radius of $\LR_{A,\phe}$, we obtain an upper bound of the essential spectral radius of $\LR_{M}$ (see the key proposition \ref{prop:conv_quasicomp}). The lower bounded is yielded by the existence of eigenvalues in open ball with the radius $\theta \lam$. Therefore the above (II) is obtained. Other results are mainly due to techniques of transfer operators in \cite{ADU, Baladi, Demers_, Lin, MU, Sar01_2, T2009}.
\smallskip
\indent
The motivation of our study is to extend the result of perturbed Markov systems with holes in the finite state case \cite{T2020,T2019} to the infinite state case. Precisely, we consider a perturbed system $(X_{A},\phe)$ with perturbed functions $\phe$ defined on $X_{A}$ satisfying that the perturbed system has a unique Gibbs measure $\mu(\e,\cd)$ of the potential $\ph(\e,\cd)$ for each $\e>0$ and on the other hand the unperturbed system $(X_{M},\ph)$ possesses several shift-invariant probability measures $(\mu_{n})_{n\geq 1}$. 
In previous work \cite{T2020,T2019} in the finite state case, we gave a necessary and sufficient condition for convergence of $\mu(\e,\cd)$ and showed that the limit measure is a convex combination of some Gibbs measures. In the infinite state case, 
however, the assumption that the unperturbed system possesses Gibbs measures imposes a restriction on the transition matrix $M$. To avoid restriction in $M$, we need an additional condition for $\ph$ and this is the summability (\ref{eq:sum}). We shall apply our results to infinite-state Markov systems with holes in a future work. In this paper, we demonstrate convergence of the topological pressure and the Gibbs measure of perturbed potential in our open system setting (Section \ref{sec:conv_pres_Gibbs}). In other examples, we apply our results to iterated function systems endowed with strongly connected multigraph (Section \ref{sec:ex2}), and treat locally constant potentials in TMS (Section \ref{sec:ex3}).
\smallskip
\indent
We recall in Section \ref{sec:Shiftsp} the notion and the results of topological Markov shifts. We state in Section \ref{sec:Ruelleop} the fundamental results of Ruelle operators and the derivation of thermodynamic features. The main results are described in Section \ref{sec:main}.
In Section \ref{sec:aux}, we give some auxiliary propositions which need to show our main results. 
The proofs of the main results are devoted in Section \ref{sec:proof}. We treat in Section \ref{sec:app_ex} various applications and examples. 
Finally, we provide a method of extending potentials defined on $X_{M}$ to on $X_{A}$ in Appendix \ref{sec:Ext_poten} and give a brief review of perturbed eigenvectors in Appendix \ref{sec:per_evec}.
\medskip
\\
{\it Acknowledgment.}\ 
The author thanks Nakano Yushi in Tokai University and Shintaro Suzuki in Tokyo Gakugei University for valuable comments and advice. This study was supported by JSPS KAKENHI Grant Number 20K03636.
\section{Preliminaries}\label{sec:pre}
\subsection{Topological Markov shifts}\label{sec:Shiftsp}
We recall the notion of topological Markov shifts. 
Let $A$ be a zero-one matrix indexed by countable set $S$. Fix a nonempty subset $S_{0}\subset S$ and an $S_{0}\times S_{0}$ zero-one matrix $M=(M(ij))$ with $M(ij)\leq A(ij)$. Consider the set
\alil
{
\textstyle{X_{M}:=\{\om\in \prod_{n=0}^{\infty}S_{0}\,:\,M(\om_{n}\om_{n+1})=1\text{ for all }n\geq 0\}}.\label{eq:XM=...}
}
Namely, $X_{M}$ is a subsystem of $X_{A}$ with the state space $S_{0}$ and the transition matrix $M$. By the sake of convenience, if $i\in S\setminus S_{0}$ or $j\in S\setminus S_{0}$, then $M(ij)$ is referred as zero.
\smallskip
\par
A word $w=w_{1}w_{2}\dots w_{n}\in S^{n}$ is {\it $M$-admissible} if $M(w_{i}w_{i+1})=1$ for all $1\leq i<n$. We write $W_{n}(M)=\{w\in S^{n}\,:\,w \text{ is }M\text{-admissible}\}$.
For $a,b\in S$, we write $a\overset{n}{\to} b$ if there exist an integer $n\geq 1$ and $w\in W_{n+1}(M)$ such that $w_{1}=a$ and $w_{n+1}=b$.
The matrix $M$ is said to be {\it irreducible} if for any $a,b\in S$, $a\overset{n}{\to} b$ for some $n\geq 1$.
We say that $M$ is {\it weakly primitive} if for any $a,b\in S$, there exists $N_{ab}\geq 1$ such that $a\overset{n}{\to} b$ for any $n\geq N_{ab}$, and $M$ is {\it primitive} if $\sup_{a,b\in S}N_{ab}<\infty$.
The matrix $M$ is said to be {\it finitely irreducible} if there exists a finite subset $F$ of $\bigcup_{n=1}^{\infty}S^{n}$ such that for any $a,b\in S$, there is $w\in F$ so that $awb$ is $M$-admissible. The matrix $M$ is called {\it finitely primitive} if there exist an integer $N\geq 1$ and a finite subset $F$ of $S^{N}$ such that for any $a,b\in S$, $awb$ is $M$-admissible. 
The matrix $M$ has the {\it big images and pre-images} (BIP) property if there is a finite subset $F\subset S$ such that for any $a\in S$, there exist $b,c\in F$ such that $M(ba)=M(ac)=1$.
Note that $M$ is finitely irreducible (finitely primitive) if and only if $M$ is irreducible (weakly primitive) and has the BIP property.
\smallskip
\par
Next we introduce tools for non-irreducible transition matrix. For $a,b\in S$, we write $a\leftrightarrow b$ when either $a\overset{m}{\to}b$, $b\overset{n}{\to}a$ for some $m,n\geq 1$ or $a=b$. Consider the quotient space $S/\!\!\leftrightarrow$. For $S_{1},S_{2}\in S/\!\!\leftrightarrow$, $S_{1}\preceq S_{2}$ denotes either $S_{1}=S_{2}$ or there exist $a\in S_{1}$ and $b\in S_{2}$ such that $a\overset{n}{\to} b$ for some $n\geq 1$. Then the relation $\preceq$ is a semi-order on $S/\!\!\leftrightarrow$. For $T\in S/\!\!\leftrightarrow$, denoted by $M(T)$ the irreducible submatrix of $M$ indexed by $T$. Therefore we obtain countable many transitive components $X_{M(T)}$ of $X$ $(T\in S/\!\!\leftrightarrow)$. 
\subsection{Ruelle transfer operators and thermodynamic formalism}\label{sec:Ruelleop}
We will recall some facts of Ruelle transfer operators which were manly introduced by \cite{MU}.
Let $X$ be a topological Markov shift with countable state space $S$ and transition matrix $A$.
For real-valued function $\ph$ on $X$, the {\it topological pressure} $P(\ph)$ of $\ph$ is given by
\alil
{
P(\ph):=\lim_{n\to \infty}\frac{1}{n}\log \sum_{w\in S^{n}\,:\,[w]\neq \emptyset}\exp(\sup_{\om\in [w]}S_{n}\ph(\om))\label{eq:toppres}
}
formally, where we put $S_{n}\ph(\om):=\sum_{i=0}^{n-1}\ph(\si^{i}\om)$. It is known in \cite[Lemma 2.1.2]{MU} that if $\ph\,:\,X\to \R$ is summable, then the number $P(\ph)$ exists in $[-\infty,+\infty)$.
\smallskip
\par
A $\si$-invariant Borel probability measure $\mu$ on $X$ is called a {\it Gibbs measure} of a potential $\ph\,:\,X\to \R$ if there exist constants $c\geq 1$ and $P\in \R$ such that for any $\om\in X$ and $n\geq 1$
\ali
{
c^{-1}\leq \frac{\mu([\om_{0}\om_{1}\dots \om_{n-1}])}{\exp(-nP+S_{n}\ph(\om))}\leq c.
}
\thms
[{\cite{MU,Sar03}}]\label{th:ex_Gibbs}
Let $X$ be a topological Markov shifts whose transition matrix is irreducible. Assume that $\ph\in F^{1}(X,\R)$ is summable. Then $\ph$ possesses an (invariant) Gibbs measure if and only if $A$ is finitely irreducible.
\thme
\pros
Due to \cite[Theorem 2.2.6 and Corollary 2.7.4]{MU}. See also \cite[Theorem 1]{Sar03}.
\proe
Theorem \ref{th:ex_Gibbs} tells that the existence of Gibbs measure imposes a restriction in the transition matrix.
\smallskip
\par
Finally we recall spectral properties of Ruelle operators $\LR_{A}=\LR_{A,\ph}$ of $\ph$. 
For $c>0$, let
\alil
{
\Lambda^{k}_{c}=\Lambda^{k}_{c}(X):=\{f\in C(X)\,:\,0\leq f,\ f(\om)\leq e^{cd_{\theta}(\om,\up)}f(\up) \text{ if }\om\in [\nu_{0}\cdots\nu_{k-1}]\}.\label{eq:Lambda_c=...}
}
The properties of $\Lambda^{k}_{c}$ will be given in Proposition \ref{prop:prop_Lamc0}.
The following is a special version of Ruelle-Perron-Frobenius Theorem for $\LR_{A}$.
\thms
[{\cite{Sar03,Sar01_2}}]\label{th:exGibbs}
Let $X$ be a topological Markov shift whose transition matrix $A$ is finitely irreducible. Assume that $\ph\in F^{1}(X,\R)$ is summable. Let $k\geq 1$ be an integer and $c\geq \theta [\ph]_{k+1}/(1-\theta)$. Then there exists a unique triplet $(\lam,h,\nu)\in \R\times \Lambda^{k}_{c}\times C_{b}(X)^{*}$ such that $\lam$ is equal to the spectral radius of $\LR_{A}\,:\,F^{1}_{b}(X)\to F^{1}_{b}(X)$, $\LR_{A}h=\lam h$, $\LR_{A}^{*}\nu=\lam\nu$ and $\nu(1)=\nu(h)=1$.
In particular, $h$ is bounded uniformly away from zero and infinity, $\lam$ equals $\exp(P(\ph))$, and $\mu=h\nu$ becomes the Gibbs measure of the potential $\ph$.
\thme
\pros
The existence of such a triplet $(\lam,h,\nu)$ is also referred a positive recurrent in \cite{Sar01_2}. Since $\ph$ has a Gibbs measure from Theorem \ref{th:ex_Gibbs}, positive recurrence of $\ph$ is guaranteed (see the proof of \cite[Theorem 1]{Sar03}). It is remain to check that $h\in \Lambda_{c}^{k}$ with $c\geq \theta [\ph]_{k+1}/(1-\theta)$. The function $h$ satisfies the inequality $\var_{k}(\log h)\leq \sum_{n=l+1}^{\infty}\var_{n}\ph$ for each $l\geq 1$ (see \cite[Corollary 2]{Sar03} or \cite[Proposition 3.4]{Sar09}). Therefore for $\om,\up$ with $\om_{0}\cdots\om_{l-1}=\up_{0}\cdots\up_{l-1}$, $\om_{l}\neq \up_{l}$ and $l\geq k$, we have $|\log h(\om)-\log h(\up)|\leq \sum_{n=l+1}^{\infty}[\ph]_{k+1}\theta^{n}=(([\ph]_{k+1}\theta)/(1-\theta)) d_{\theta}(\om,\up)\leq cd_{\theta}(\om,\up)$. Hence $h\in \Lambda_{c}^{k}$. 
\proe
\section{Main results}\label{sec:main}
The following is one of our main results.
\thms
\label{th:ex_efunc_geneRuelleop_M}
Assume that the conditions (A.1)-(A.3) are satisfied. Put $c_{\adl{Lamc}}:=[\ph]_{k+1}\theta/(1-\theta)$. Then we have the following:
\ite
{
\item There exists a triplet $(\lam,g,\nu)$ such that (i) $\lam$ is equal to $\exp(P(\ph|_{X_{M}}))$ and is the positive spectral radius of the operator $\LR_{M}\,:\,F^{k}_{b}(X)\to F^{k}_{b}(X)$ and of $\LR_{M}\,:\,C_{b}(X)\to C_{b}(X)$ both, (ii) $g\in \Lambda^{k}_{c_{\adr{Lamc}}}(X)$ with $\|g\|_{\infty}=1$ and $\LR_{M}g=\lam g$, and (iii) $\nu$ is a Borel probability measure on $X$ supported on $X_{M}$ and $\LR_{M}^{*}\nu=\lam \nu$. 
\item If $X_{A}$ is not reduced to finite single orbit, then the essential spectral radius $r_{\text{ess}}(\LR_{M})$ of the operator $\LR_{M}\,:\,F^{k}_{b}(X)\to F^{k}_{b}(X)$ satisfies $\lam\theta\leq r_{\text{ess}}(\LR_{M})<\lam$. In particular, for any $r\in [0,\lam\theta)$ except for at most a countable number, $p\in \C$ with $|p|=r$ is an eigenvalue with infinite multiplicity. If $X_{A}$ is a finite single orbit, then $r_{\text{ess}}(\LR_{M})=0$.
}
\thme
For the corollary of this theorem, we introduce the following condition which is stronger than the condition (A.1):
\ite
{
\item[(A.1)${}^\p$] The matrix $A$ is finitely irreducible.
}
\cors\label{cor:ex_efunc_geneRuelleop_M_fi}
Assume that the conditions (A.1)${}^\p$, (A.2) and (A.3) are satisfied. Then if $X_{A}$ is not reduced to finite single orbit, then the essential spectral radius $\LR_{M}\,:\,F^{k}_{b}(X)\to F^{k}_{b}(X)$ is equal to $\lam\theta$.
\core
The proofs of Theorem \ref{th:ex_efunc_geneRuelleop_M} and Corollary \ref{cor:ex_efunc_geneRuelleop_M_fi} will be stated in Section \ref{sec:proof_th:ex_efunc_geneRuelleop_M}.
\rem
{
\ite
{
\item Assume that the conditions (A.1)-(A.3) are satisfied. The TMS $X$ can be recoded in a natural way to a TMS with the state $S^{*}=\{w\in S^{k+1}\,:\, M\text{-admissible}\}$ and with the transition matrix $(A^{*}(ww^\p))$ as $A^{*}(ww^\p)=1$ iff $w_{2}\cdots w_{k}=w_{1}^\p\cdots w_{k-1}^\p$. Then $\ph$ is reduced to a potential $\ph^{*}\,:\,X_{A^{*}}\to \R$ satisfying $[\ph^{*}]_{1}<\infty$. However $\ph^{*}$ may not be summable.
\item Assume that (A.1) and (A.3) are satisfied. Then (A.2) holds if and only if the spectral radius of $\LR_{M}$ is positive (Proposition \ref{prop:positive_eval}).
Moreover, if the state space $S$ is finite, then the condition (A.2) is satisfied if and only if $X_{M}$ is non-empty.
}
}
To state our second result, assume that $M$ is irreducible. In fact, since any nonnegative matrix is decomposed into topological transitive components, it is important to know spectrum properties of $\LR_{M}$ when $M$ is irreducible. 
Denoted by $p$ the period of $M$. Since $M$ is irreducible, there exists a decomposition $S_{0}=S_{0,0}\cup\cdots\cup S_{0,p-1}$ such that
\alil
{
\textstyle{\Si_{i}:=\bigcup_{s\in S_{0,i}}[s]},\qquad
X_{M}^{i}:=X_{M}\cap \Si_{i},\qquad \si(X_{M}^{i})=X_{M}^{i+1\text{ mod } p}\qquad (0\leq i<p )\label{eq:Sii=}
}
and each $(X_{M}^{i},\si^{p})$ is topologically mixing.
Then we have the following:
\thms
\label{th:ex_efunc_geneRuelleop_M_irre}
Assume that the conditions (A.1)-(A.3) are satisfied and the matrix $M$ is irreducible. Let $p\geq 1$ be the period of $M$ and put $\kappa=\exp(2\pi\sqrt{-1}/p)$. Then the operator $\LR_{M}\,:\,F^{k}_{b}(X)\to F^{k}_{b}(X)$ obtain the spectral decomposition
\alil
{
\textstyle{\LR_{M}=\sum_{i=0}^{p-1}\lam_{i}\PR_{i}+\RR}\label{eq:LMph=sumi=0p-1...}
}
satisfying that each $\lam_{i}=\lam\kappa^{i}$ is a simple eigenvalue of $\LR_{M}$, $\PR_{i}$ is a projection and has the form $\PR_{i}=h_{i}\otimes \nu_{i}$ with the eigenfunction $h_{i}=\sum_{j=0}^{p-1}\kappa^{-ji}h\chi_{\Si_{j}}$ and the eigenvector $\nu_{i}=\sum_{j=0}^{p-1}\kappa^{ji}\nu|_{\Si_{j}}$, and the spectral radius of $\RR$ is less than $\lam$, where $h=g/\nu(g)$ and $(\lam,g,\nu)$ appears in Theorem \ref{th:ex_efunc_geneRuelleop_M}(1). Here $\chi_{\Si_{j}}$ means the indicator of the set $\Si_{j}$.
\thme
\cors
\label{cor:ex_efunc_geneRuelleop_M_irre}
Assume that the conditions (A.1)-(A.3) are satisfied and the set $\{S_{1}\in S/\!\!\leftrightarrow\,:\,P(\ph|_{X_{M(S_{1})}})=P(\ph|_{X_{M}})\}$ is only one element $S_{1}$. Let $p$ be the period of $M(S_{1})$. 
Then the operator $\LR_{M}\,:\,F^{k}_{b}(X)\to F^{k}_{b}(X)$ has the spectral decomposition as well as (\ref{eq:LMph=sumi=0p-1...}). In particular, $\supp\, h_{i}=\bigcup\{\bigcup_{a\in T}[a]\,:\,T\preceq S_{1}\}$ and $\supp\,\nu_{i}=X_{M}\cap \{\bigcup_{a\in T}[a]\,:\,S_{1}\preceq T\}$ are satisfied.
\core
We prove Theorem \ref{th:ex_efunc_geneRuelleop_M_irre} and Corollary \ref{cor:ex_efunc_geneRuelleop_M_irre}
in Section \ref{sec:proof_th:ex_efunc_geneRuelleop_M_irre} and Section \ref{sec:proof_cor:ex_efunc_geneRuelleop_M_irre}, respectively.
\rem
{\label{rem:SPR}
Let $X_{A}$ be a topological Markov shift whose shift is topological mixing and $\ph\,:\,X_{A}\to \R$ a weak $d_{\theta}$-Lipschitz continuous and summable, and we put $M=A$ (i.e. the cases $X_{A}=X_{M}$, $p=1$ and $k=1$ in Theorem \ref{th:ex_efunc_geneRuelleop_M_irre}). Then the result of Theorem \ref{th:ex_efunc_geneRuelleop_M_irre} tells us that $\ph$ is strongly positive recurrent (SPR) in the sense of the notion in \cite{Cyr_Sarig}. Precisely, if we define a Banach space $F=\{f\in C(X_{M}):\|f\|_{F}<\infty\}$ endowed with the norm $\|\cd\|_{F}$:
\ali
{
\|f\|_{F}:=&\sup_{b\in S}\{(\sup_{[b]}h)^{-1}(\sup_{[b]}|f|+\sup\{|f(\om)-f(\up)|/\theta^{s_{a}(\om,\up)}\,:\,\om,\up\in [b],\ \om\neq \up\})\}\\
s_{a}(\om,\up):=&\sharp\{0\leq i\leq t(\om,\up)-1\,:\,x_{i}=y_{i}=a\}\quad \text{with}\quad t(\om,\up):=\min\{n\,:\,\om_{n}\neq \up_{n}\},
}
where fix $a\in S$, then the Ruelle operator $\LR_{A}\,:\,F\to F$ has a spectral gap property.
In the case $\inf_{X} h=0$, any bounded function $f\in F^{1}_{b}(X)$ with $\inf_{X} |f|>0$ is not in $F$. Thus $F^{1}_{b}(X)\nsubseteq F$ in general. Theorem \ref{th:ex_efunc_geneRuelleop_M_irre} states that the Ruelle operator $\LR_{A}$ also has a spectral gap property on $F^{1}_{b}(X)$. 
}
Finally, we give the escape rate from the open system $\si|_{\Si}$ under the conditions (A.1)-(A.3). Recall the set $\Si=\bigcup_{ij\in S^{2}\,:\,M(ij)=1}[ij]$ and put $\Si^{n}=\bigcap_{i=0}^{n}\si^{-i}\Si$. By Theorem \ref{th:ex_efunc_geneRuelleop_M_irre} replacing $M$ by $A$, there exists a unique triplet $(\lam_{A},h_{A},\nu_{A})$ such that $\lam_{A}$ is spectral radius of $\LR_{A}$ acting on $F^{k}_{b}(X)$ and is equal to $\exp(P(\ph))$, $\LR_{A}h_{A}=\lam_{A}h_{A}$, $\LR_{A}^{*}\nu_{A}=\lam_{A}\nu_{A}$ and $\nu_{A}(1)=\nu_{A}(h_{A})=1$. Put $\mu_{A}=h_{A}\nu_{A}$. Then $\mu_{A}$ is a $\si$-invariant Borel probability measure and $\supp\,\mu_{A}=X$ (see Proposition \ref{prop:supp_nu}).
Then we obtain the following:
\thms
\label{th:escaperate}
Assume that the conditions (A.1)-(A.3) are satisfied. Then we have the limit $\lim_{n\to \infty}(1/n)\log \mu_{A}(\Si^{n})=P(\ph|_{X_{M}})-P(\ph)$.
\thme
We will give this proof in Section \ref{sec:escaperate}. When $A\neq M$, $P(\ph|_{X_{M}})$ is strictly less than $P(\ph)$. Therefore the escape rate of the open system $\si|_{\Si}$ is exponential.
\section{Auxiliary propositions}\label{sec:aux}
In this section, we enumerate and show some auxiliary propositions which are useful to prove our main results. We use the notation defined in Section \ref{sec:intro} and in Section \ref{sec:pre}.
\prop
{\label{prop:prop_Lamc0}
Assume that the conditions (A.1)-(A.3) are satisfied. Then
\ite
{
\item if $w\in \bigcup_{i=1}^{k}S^{i}$ then $\chi_{[w]}\in \Lambda^{k}_{c}$ for all $c\geq 0$, where $\chi_{[w]}$ is the indicator of the set $[w]$;
\item if $f,g\in \Lambda^{k}_{c}$ and $\alpha,\beta\geq 0$, then $\alpha f+\beta g\in \Lambda^{k}_{c}$;
\item if $c\geq c_{\adr{Lamc}}$ and $f\in \Lambda^{k}_{c}$ with $\LR_{M}f\in C(X)$, then $\LR_{M}f \in \Lambda^{k}_{c}$.
}
}
\pros
The assertions (1)-(2) immediately follow by direct checking. We will see the validity of (3). For $f\in \Lambda^{k}_{c}$ and $d(\om,\up)\leq \theta^{k}$
\ali
{
\LR_{M} f(\om)
\leq&\sum_{a\,:\,M(a\up_{0})=1}e^{\ph(a\cd\up)+[\ph]_{k+1}d_{\theta}(a\cd \om,a\cd\up)}f(a\cd\up)e^{cd_{\theta}(a\cd\om,a\cd\up)}=\LR_{M}f(\up) e^{\theta([\ph]_{k+1}+c)d_{\theta}(\om,\up)}.
}
Here we note $\theta([\ph]_{k+1}+c)\leq c$ if and only if $c_{\adr{Lamc}}\leq c$. Hence $\LR_{M} f\in \Lambda^{k}_{c}$.
\proe
\noindent
Denoted by $r(\LR)$ the spectral radius of a bounded linear operator $\LR$ acting on $F^{k}_{b}(X)$, and by $r_{C}(\LR)$ the spectral radius of $\LR$ acting on $C_{b}(X)$.
\prop
{\label{prop:positive_eval}
Assume that the conditions (A.1) and (A.3) are satisfied. Assume also that $M$ is an $S_{0}\times S_{0}$ ($S_{0}\subset S$) zero-one matrix such that $M(ij)\leq A(ij)$. Then $r(\LR_{M})=r_{C}(\LR_{M})$. Moreover $r(\LR_{M})$ is positive if and only if (A.2) is satisfied.
}
\pros
First we show $r(\LR_{M})=r_{C}(\LR_{M})$. By $\|\LR_{M}^{n}\|_{\infty}=\|\LR_{M}^{n}1\|_{\infty}\leq \|\LR_{M}^{n}1\|_{k}$, we have $r_{C}(\LR_{M})\leq r(\LR_{M})$. To see the inverse inequality, let $n\geq k$, $f\in F^{k}_{b}(X)$ and $\om,\up\in X$ with $d_{\theta}(\om,\up)\leq \theta^{k}$. We note the basic inequality
\alil
{
|\LR_{M}^{n}f(\om)-\LR_{M}^{n}f(\up)|
\leq&\sum_{w\,:\,w\cd\om_{0}\in W_{n+1}(M)}\{e^{S_{n}\ph(w\cd\om)}|f(w\cd\om)-f(w\cd\up)|+\label{eq:basine}\\
&\qqqquad e^{S_{n}\ph(w\cd\up)}|e^{S_{n}\ph(w\cd\om)-S_{n}\ph(w\cd\up)}-1||f(w\cd\up)|\}\nonumber\\
\leq&\|\LR_{M}^{n}1\|_{\infty}([f]_{k}\theta^{n}+c_{\adr{basinq}}\|f\|_{\infty})d_{\theta}(\om,\up)\nonumber
}
with $c_{\adl{basinq}}=e^{\theta^{k}c_{\adr{Lamc}}}c_{\adr{Lamc}}$.
Therefore $[\LR_{M}^{n}f]_{k}\leq (1+c_{\adr{basinq}})\|\LR_{M}^{n}1\|_{\infty}\|f\|_{k}$. Consequently, we obtain $\|\LR_{M}^{n}\|_{k}\leq (2+c_{\adr{basinq}})\|\LR_{M}^{n}1\|_{\infty}$ and thus $r(\LR_{M})\leq r_{C}(\LR_{M})$.
\smallskip
\par
Next we prove the second assertion. Assume that $X_{M}$ has a period point $\om$ with $\si^{m}\om=\om$.
For any integer $l\geq 1$, we have
\alil
{
\LR_{M}^{lm}1(\om)=&\sum_{w\in (S_{0})^{lm}\,:\,w\cd \om_{0}\in W_{lm+1}(M)}e^{S_{lm}\ph(w\cd \om)}\geq e^{S_{lm}\ph(\om)}=e^{lS_{m}\ph(\om)}.\label{eq:LMph^km1=}
}
Thus $r_{c}(\LR_{M})\geq e^{S_{m}\ph(\om)/m}>0$. Hence we obtain the assertion.

Conversely, we assume that $X_{M}$ has no period point. Choose any $\eta>0$. Since $\ph$ is summable, there exists a subset $S_{1}\subset S$ such that $S_{2}:=S\setminus S_{1}$ is finite and $\sum_{s\in S_{1}}e^{\sup_{[s]}\ph}<\eta$. Since any finite $M$-admissible $w=w_{1}\cdots w_{n}$ satisfies that $w_{1},\dots, w_{n}$ are different each other, we have $\sharp\{i\,:\,w_{i}\in S_{2}\}\leq m:=\sharp S_{2}$. We see for $n>m$
\ali
{
\LR_{M}^{n}1(\om)
\leq&\sum_{w_{1},\dots,w_{n}\in S\,:\,\text{mutually different}}e^{\sum_{i=1}^{n}\sup_{[w_{i}]}\ph}\\
\leq&\sum_{w_{1},\dots,w_{n-m}\in S_{1}\atop{w_{n-m+1},\dots,w_{n}\in S}}e^{\sum_{i=1}^{n}\sup_{[w_{i}]}\ph}
\leq \Big(\sum_{s\in S_{1}}e^{\sup_{[s]}\ph}\Big)^{n-m}\Big(\sum_{s\in S}e^{\sup_{[s]}\ph}\Big)^{m}\leq c\eta^{n-m}
}
for a constant $c>0$. Thus we obtain $r(\LR_{M})\leq \eta$.
Hence $r_{c}(\LR_{M})=0$.
\proe
Assume that the conditions (A.1)-(A.3) are satisfied. We also introduce the following conditions for perturbed potentials $\phe$ with a small parameter $\e\in (0,1)$:
\ite
{
\item[(B.1)] $\ph(\e,\cd)\in F^{1}(X,\R)$ for any $\e>0$ and $c_{\adl{udphe}}:=\sup_{\e>0}[\phe]_{k+1}$ is finite.
\item[(B.2)] $c_{\adl{unisum}}:=\sum_{s\in S}\exp(\sup_{\e>0}\sup_{\om\in [s]}\ph(\e,\om))<\infty$.
\item[(B.3)] Take $\psi\,:\,X\to \R$ so that $\psi(\om)=\exp(\ph(\om))$ if $M(\om_{0}\om_{1})=1$ and $\psi(\om)=0$ otherwise. Then $\sup_{\om\in [a]}|e^{\ph(\e,\om)}-\psi(\om)|\to 0$ as $\e\to 0$ for each $a\in S$.
}
For convenience, let $\lam(\e)=\exp(P(\phe))$, $c_{\adl{Lc}}=c_{\adr{udphe}}\theta/(1-\theta)$ and $c_{\adl{udphe2}}=\max(c_{\adr{Lc}},[\ph]_{k+1}\theta/(1-\theta))$.
Take $\up\in X_{M}$ so that $\si^{m}\up=\up$ and put
$\lam_{*}=(\prod_{i=0}^{m-1}\psi(\si^{i}\up))^{1/m}$.
Then we notice $\liminf_{\e\to 0}\lam(\e)\geq \lam_{*}$ by (\ref{eq:LMph^km1=}) and (B.3).
\prop
{\label{prop:conv_Ruelleop_LB}
Assume that the conditions $(A.1)$-$(A.3)$ and $(B.1)$-$(B.3)$ are satisfied.
Then $\|\LR_{A,\phe}-\LR_{M}\|_{\infty}\to 0$ as $\e\to 0$.
}
\pros
Put $N=\bigcup_{ij\,:\,A(ij)=1,M(ij)=0}[ij]$. We write $S=\{s_{1},s_{2},\dots, s_{d}\}$ with $d\leq \infty$. For $n=1,2,\dots$, $\e> 0$, $\om\in X$ and $f\in C_{b}(X)$, we define
\ali
{
H_{n}(\e,\om,f)=&e^{\ph(\e,s_{n}\cd\om)}f(s_{n}\cd\om),\quad
H_{n}(0,\om,f)=\psi(s_{n}\cd\om)(1-\chi_{N}(s_{n}\cd\om))f(s_{n}\cd\om)
}
if $A(s_{n}\om_{0})=1$, and $H_{n}(\e,\om,f)=H_{n}(0,\om,f)=0$ if $A(s_{n}\om_{0})=0$.
Then we see
$\LR_{A,\phe}f(\om)=\sum_{n=1}^{\infty}H_{n}(\e,\om,f)$.
Put $a_{n}=\exp(\sup_{\e>0}\sup_{\om\in [s_{n}]}\ph(\e,\om))+\sup_{\om\in [s_{n}]}\psi(\om)$. Note $\sum_{n=1}^{\infty}a_{n}<\infty$ by the conditions $(B.2)$ and the summability of $\ph$. We obtain
$|H_{n}(\e,\om,f)|\leq a_{n}\|f\|_{\infty}$ for any $\e\geq 0$, $\om$ and $f$. Moreover
\ali
{
\textstyle |H_{n}(\e,\om,f)-H_{n}(0,\om,f)|\leq \sup_{\om\in [s_{n}]}|e^{\ph(\e,\om)}-\psi(\om)(1-\chi_{N}(\om))|\|f\|_{\infty}.
}
Choose any $\eta>0$. Then there exists $n_{0}\geq 1$ such that
$\sum_{n=n_{0}+1}^{\infty}a_{n}<\eta$.
Furthermore, the condition (B.3) implies that there exists $\e_{1}>0$ such that for any $0<\e<\e_{1}$ and for any $a\in \{s_{1},s_{2},\dots,s_{n_{0}}\}$, 
$\sup_{\om\in [a]}|e^{\ph(\e,\om)}-\psi(\om)(1-\chi_{N}(\om))|<\eta/n_{0}$. 
Therefore, for any $f\in C_{b}(X)$ with $\|f\|_{\infty}\leq 1$ and $0<\e<\e_{1}$
\ali
{
&\|(\LR_{A,\phe}-\LR_{M})f\|_{\infty}=\sup_{\om\in X}|\sum_{n\geq 1}H_{n}(\e,\om,f)-\sum_{n\geq 1}H_{n}(0,\om,f)|\\
\leq&\sup_{\om\in X}\sum_{1\leq n\leq n_{0}}|H_{n}(\e,\om,f)-H_{n}(0,\om,f)|+\sup_{\om\in X}\sum_{n\geq n_{0}+1}(|H_{n}(\e,\om,f)|+|H_{n}(0,\om,f)|)
\leq 3\eta.
}
This implies that $\LR_{A,\phe}$ converges $\LR_{M}$ with respect to $\|\cd\|_{\infty}$.
\proe
\prop
{\label{prop:conv_evec_LB}
Assume that the conditions $(A.1)$-$(A.3)$ and $(B.1)$-$(B.3)$ are satisfied. Assume also that for each $\e>0$, there exists a Borel probability measure $\nu(\e,\cd)$ on $X$ such that $\LR_{A,\phe}^{*}\nu(\e,\cd)=\lam(\e)\nu(\e,\cd)$. Then $\{\nu(\e,\cd)\}$ has a converging subsequence in the sense of weakly topology. Moreover, any limit point $\nu$ of $\nu(\e,\cd)$ as $\e\to 0$ is a Borel probability measure on $X$ and an eigenvector of $\lam$ of the dual $\LR_{M}^{*}$ of $\LR_{M}\,:\,C_{b}(X)\to C_{b}(X)$, where $\lam$ is a limit point of $\lam(\e)$.
}
\pros
We write $S=\{1,2,\dots, d\}$ with $d\leq \infty$.
First we will use a technique of Prohorov's theorem for the sequence $\{\nu(\e,\cd)\}$, i.e. we will show the tightness of this. 
For $n\geq 1$, let $\pi_{n}\,:\,X\to S$ be $\pi_{n}(\om)=\om_{n}$. In this case, $\pi_{n}$ is continuous. For $s\in S$ we have
\ali
{
\lam(\e)^{n}\nu(\e,\pi_{n}^{-1}(s))
=&\int_{X}\sum_{w\in S^{n}\,:\,w\cd\om_{0}\in W_{n+1}(A)}e^{S_{n}\ph(\e,w\cd\om)}\chi_{\pi_{n}^{-1}(s)}(w\cd\om)\,d\nu(\e,\om)\leq c_{\adr{unisum}}^{n-1}e^{\sup_{[s]}\ph(\e,\cd)},
}
where $c_{\adr{unisum}}$ appears in (B.2). Choose any $\eta>0$. By using $\lam(\e)\geq \lam_{*}/2$ for any small $\e>0$,
\ali
{
\textstyle\nu(\e,\pi_{n}^{-1}([s+1,\infty)))\leq (\lam_{*}/2)^{-n}c_{\adr{unisum}}^{n-1}\sum_{j>s}e^{\sup_{\e>0}\sup_{[j]}\ph(\e,\cd)}.
}
Therefore $\nu(\e,\pi_{n}^{-1}([s(n)+1,\infty)))\leq \eta/2^{n}$ for some $s(n)\geq 1$.
Thus
\ali
{
\textstyle\nu(\e,\bigcap_{n=1}^{\infty}\pi_{n}^{-1}[1,s(n)])=1-\nu(\e,\bigcup_{n=1}^{\infty}\pi_{n}^{-1}[s(n)+1,\infty))\geq 1-\eta.
}
Since $\bigcap_{n=1}^{\infty}\pi_{n}^{-1}[1,s(n)]$ is compact, the sequence $\{\nu(\e,\cd)\}$ is tight. Prohorov's theorem implies that there exist a subsequence $\{\nu(\e_{n},\cd)\}$ and a Borel probability measure $\nu$ on $X$ such that $\nu(\e_{n},\cd)$ converges to $\nu$ weakly, namely $\nu(\e_{n},f)\to \nu(f)$ as $n\to \infty$ for each $f\in C_{b}(X)$.
Since $\{\lam(\e)\}$ is bounded, we may assume convergence $\lam_{\e_{n}}\to \lam$. We have
\ali
{
&|\lam\nu(f)-\LR_{M}^{*}\nu(f)|\\ 
\leq&|\lam\nu(f)-\lam(\e)\nu(\e,f)|+|\nu(\e,\LR_{A,\phe}f)-\nu(\e,\LR_{M}f)|+|\nu(\e,\LR_{M}f)-\nu(\LR_{M}f)|\\
\leq&|\lam\nu(f)-\lam(\e)\nu(\e,f)|+\|\LR_{A,\phe}-\LR_{M}\|_{\infty}\|f\|_{\infty}+|(\nu(\e,\cd)-\nu)(\LR_{M}f)|
\to 0
}
as $\e\to 0$ running through $\{\e_{n}\}$. Hence the proof is complete.
\proe
\prop
{\label{prop:supp_nu}
Assume that the conditions $(A.1)$-$(A.3)$ are satisfied and $M$ is irreducible. If $\nu\neq 0$ is a finite Borel measure with $\LR_{M}^{*}\nu=\lam\nu$ and $\lam>0$, then $\supp\, \nu=X_{M}$.
}
\pros
Choose any $\om\in X\setminus X_{M}$. There is $m\geq 1$ so that $M(\om_{m-1}\om_{m})=0$ and therefore $\nu([\om_{0}\cdots \om_{m}])=\lam^{-m-1}\nu(\LR_{M}^{m+1}\chi_{[\om_{0}\cdots \om_{m}]})=0$. Thus $\supp\ \nu\subset X_{M}$. To check the converse, we show $\nu(O)>0$ for any open set $O$ with $\om\in O$ for each $\om\in X_{M}$. Let $\tau\in S_{0}^{k}$ with $\nu([\tau])>0$. Since $M$ is irreducible, for any $n\geq 0$ there exists $w\in S_{0}^{m}$ such that $\om_{n}\cd w\cd \tau$ is $M$-admissible. Let $l=n+m+1$. For an element $\up\in [\tau]$, we see
\alil
{
\nu([\om_{0}\cdots \om_{n}])=&\lam^{-l}\nu(\LR_{M}^{l}\chi_{[\om_{0}\cdots \om_{n}]})
\geq \lam^{-l}e^{\theta^{k}c_{\adr{Lamc}}}e^{S_{l}\ph(\om_{0}\cdots \om_{n}\cd w\cd \up)}\nu([\tau])>0.\label{eq:hnu[...]>0}
}
Since all cylinders are open base, we get $\nu(O)>0$. Hence $\supp\ \nu=X_{M}$.
\proe
\prop
{\label{prop:conv_efunc_LB}
Assume that the conditions $(A.1)$-$(A.3)$ and $(B.1)$-$(B.3)$ are satisfied. Assume also that for each $\e>0$, there exists a nonnegative functions $g(\e,\cd)\in \Lambda^{k}_{c_{\adr{udphe2}}}$ such that $\|g(\e,\cd)\|_{\infty}=1$ and $\LR_{A,\phe}g(\e,\cd)=\lam(\e)g(\e,\cd)$. Then $\{g(\e,\cd)\}$ has a converging subsequence $\{g(\e_{n},\cd)\}$ in the sense that $g(\e_{n},\om)$ converges to a function $g(\om)$ for each point $\om\in X$. Moreover, the limit $g$ is a nonzero nonnegative function belonging in $\Lambda^{k}_{c_{\adr{udphe2}}}$ and is an eigenfunction of $\lam$ of the operator $\LR_{M}$, where $\lam$ is a limit point of $\lam(\e)$.
}
\pros
By $g(\e,\om)\leq e^{c_{\adr{Lc}}d_{\theta}(\om,\up)}g(\e,\up)$ for any $d_{\theta}(\om,\up)\leq \theta^{k}$ and by $\|g(\e,\cd)\|_{\infty}=1$, 
$\{g(\e,\cd)\}$ is equicontinuous and uniformly bounded. By Ascoli Theorem, there exist a subsequence $(\e_{n})$ and $g\in C_{b}(X)$ such that $g(\e_{n},\om)\to g(\om)$ for each $\om\in X$. From $g$ satisfies $g(\om)\leq e^{c_{\adr{Lc}}d_{\theta}(\om,\up)}g(\up)$ for $\om\in [\up_{0}\cdots\up_{k-1}]$, $g$ is in $\Lambda^{k}_{c_{\adr{Lc}}}$. We show $g\neq 0$.
There exists $S_{0}\subset S$ such that $S_{1}:=S\setminus S_{0}$ is finite and $\sum_{s\in S_{0}}e^{\sup_{[s]}\phe}<c_{\adr{unisum}}^{-k+1}(\lam_{*})^{k}/(6k)$.
For $\om\in X$
\ali
{
\frac{\lam_{*}^{k}}{2}g(\e,\om)\leq \lam(\e)^{k}g(\e,\om)=\LR_{A,\phe}^{k}g(\e,\om)&=\sum_{w\in (S_{0}\cup S_{1})^{k}\,:\,w\cd\om_{0}\in W_{k+1}(A)}e^{S_{k}\ph(\e,w\cd\om)}g(\e,w\cd\om)\\
&\leq\sum_{w\in S_{1}^{k}\,:\,w\cd\om_{0}\in W_{k+1}(A)}e^{S_{k}\ph(\e,w\cd\om)}g(\e,w\cd\om)+\frac{(\lam_{*})^{k}}{6}
}
holds for any small $\e>0$ so that $(\lam_{*})^{k}/2\leq \lam(\e)^{k}$. Choose any $\up \in X$ with $g(\e,\up)>2/3$ and $\up^{w}\in [w]$ for any $w\in W_{k}(A)\cap S_{1}^{k}$. We have
\ali
{
\frac{(\lam_{*})^{k}}{6}<\frac{(\lam_{*})^{k}}{2}g(\e,\up)-\frac{(\lam_{*})^{k}}{6}\leq e^{c_{\adr{udphe}}\theta}\sum_{w\in W_{k}(A)\cap S_{1}^{k}}\prod_{i=1}^{k}e^{\sup_{\eta>0}\sup_{[w_{i}]}\ph(\eta,\cd)}g(\e,\up^{w}).
}
Since $S_{1}^{k}$ is finite, there exist a subsequence $(\e_{n})$ and $w\in S_{1}^{k}$ such that $\inf_{n}g(\e_{n},\up^{w})>0$. this implies that any limit point $g$ of $g(\e,\cd)$ is not zero.

Take limit points $\lam$ of $\lam(\e)$ and $g$ of $g(\e,\cd)$. We show $\LR_{M}g=\lam g$. 
By virtue of Proposition \ref{prop:conv_evec_LB} replacing $M$ by $A$, there exists a pair $(\hat{\lam},\hat{\nu})$ such that $\hat{\lam}$ is positive and $\hat{\nu}$ is a Borel probability measure $\hat{\nu}$ on $X$ and satisfies
$\LR_{A}^{*}\hat{\nu}=\hat{\lam}\hat{\nu}$. 
Remark that $\hat{\nu}$ has full support on $X$ by Proposition \ref{prop:supp_nu}. We have
\ali
{
\|\lam(\e)g(\e,\cd)-\LR_{M}g\|_{L^{1}(\hat{\nu})}\leq& \|(\LR_{A,\phe}-\LR_{M})g(\e,\cd)\|_{L^{1}(\hat{\nu})}+\|\LR_{M}g(\e,\cd)-\LR_{M}g\|_{L^{1}(\hat{\nu})}\\
\leq&\|\LR_{A,\phe}-\LR_{M}\|_{\infty}+\int_{X}\LR_{A}|g(\e,\cd)-g|(\om)\,d\hat{\nu}(\om)\\
\leq& \|\LR_{A,\phe}-\LR_{M}\|_{\infty}+\hat{\lam}\|g(\e,\cd)-g\|_{L^{1}(\hat{\nu})},
}
By Lebesgue dominated convergence theorem, $g(\e_{n},\cd)$ converges to $g$ in $L^{1}(X,\hat{\nu})$ for some sequence $(\e_{n})$. 
Letting $\e=\e_{n}\to 0$, we obtain $\LR_{M}g=\lam g$ $\hat{\nu}$-a.e. 
By continuity of $g$ and $\supp\,\hat{\nu}=X$, we have $\LR_{M}g=\lam g$.
\proe
\prop
{\label{prop:supp_g}
Assume that the conditions $(A.1)$-$(A.3)$ are satisfied and $M$ is irreducible. If $g\in \Lambda^{k}_{c}$ satisfies $g\neq 0$ and $\LR_{M}g=\lam g$ with $\lam>0$, then $\supp\, g=\bigcup_{a\in S_{0}}[a]$.
}
\pros
For $\om\in X$ with $\om_{0}\notin S_{0}$, $g(\om)=\lam^{-1}\LR_{M}g(\om)=0$ by the definition of $\LR_{M}$. Therefore $\supp\, g\subset \bigcup_{s\in S_{0}}[s]$. Conversely, choose any $\om\in X$ with $\om_{0}\in S_{0}$. Take $\up\in X$ so that $w:=\up_{0}\dots\up_{k-1}\in S_{0}^{k}$ and $g(\up)>0$. 
Since $M$ is irreducible, there exists $w^\p\in S_{0}^{m-k}$ with $m>k$ such that $w\cd w^\p\cd \om_{0}$ is $M$-admissible. We see
\ali
{
g(\om)=\lam^{-m}\LR_{M}^{m}g(\om)\geq& \lam^{-m}e^{S_{m}\ph(w\cd w^\p\cd\om)}g(w\cd w^\p\cd\om)
\geq\lam^{-m}e^{S_{m}\ph(w\cd w^\p\cd\om)}e^{-c\theta^{k}}g(\up)>0,
}
where the second inequality holds from $g\in \Lambda_{c}^{k}$. Thus we obtain $\supp\ g=\bigcup_{s\in S_{0}}[s]$.
\proe
Next we show an important proposition about quasi-compactness of $\LR_{M}$.
\prop
{\label{prop:conv_quasicomp}
Assume that the conditions $(A.1){}^{\p}$, $(A.2)$ and $(A.3)$ are satisfied. We also assume that a function $\ph_{0} \in F^{1}(X,\R)$ is finite pressure and satisfies $\ph\leq \ph_{0}$ on $\{\om\in X\,:\,M(\om_{0}\om_{1})=1\}$. Let $\lam_{0}=\exp(P(\ph_{0}))$. 
Then the essential spectral radius $r_{\text{ess}}(\LR_{M})$ of the operator $\LR_{M}\,:\,F^{k}_{b}(X)\to F^{k}_{b}(X)$ satisfies $r_{\text{ess}}(\LR_{M})\leq \lam_{0}\theta$. 
}
\pros
Choose any $\eta>0$. By virtue of Theorem \ref{th:exGibbs}, there exists a unique triplet $(\lam_{0},h_{0},\nu_{0})$ such that $\LR_{A,\ph_{0}}h_{0}=\lam_{0}h_{0}$, $\LR_{A,\ph_{0}}^{*}\nu_{0}=\lam_{0}\nu_{0}$, $\nu_{0}(h_{0})=\nu_{0}(1)=1$, $\|h_{0}\|_{\infty}<\infty$ and $\|h_{0}^{-1}\|_{\infty}<\infty$ hold.
Since $A$ is finitely irreducible, $\mu_{0}:=h_{0}\nu_{0}$ is the Gibbs measure for the potential $\ph_{0}$.
Put $\QR:=\lam_{0}^{-1}\LR_{M}$.
We will show that the operator $\QR$ satisfies the Hennion's condition \cite{HH}. 
In order to this, we check the following four claims.
\cla
{\label{cla:LYineqQ}
There exist constants $c_{\adl{Qct1}},c_{\adl{Qct2}}\geq 1$ such that $\|\QR^{m}f\|_{k}\leq c_{\adr{Qct1}}\|f\|_{L^{1}(\mu_{0})}+c_{\adr{Qct2}}\|f\|_{k}\theta^{m}$ for any $f\in F^{k}_{b}(X)$ and $m\geq k$.
}
For any $f\in F^{k}_{b}(X)$, $m\geq k$, $A$-admissible word $w\in S^{m}$ and $\om,\up\in [w]$, we have
$|f(\om)-f(\up)|\leq [f]_{k}d_{\theta}(\om,\up)\leq[f]_{k}\theta^{m}$
and therefore $f(\om)\leq f(\up)+[f]_{k}\theta^{m}$.
By integrating for the measure $\mu_{0}$ on $\up\in [w]$ in both sides and by dividing by $\mu_{0}[w]$, we obtain
\ali
{
f(\om)\leq \frac{1}{\mu_{0}([w])}\int_{[w]}f(\up)\,d\mu_{0}(\up)+[f]_{k}\theta^{m}
}
for any $\om\in [w]$. 
Moreover, by the Gibbs property for $\mu_{0}$, we get $e^{S_{m}\ph_{0}(\om)-m\log \lam_{0}}\leq c_{\adr{Gibbsct1}}\mu_{0}([w])$ 
for some constant $c_{\adl{Gibbsct1}}\geq 1$. We see for $\om\in X$
\ali
{
(\QR^{m}|f|)(\om)=\lam_{0}^{-m}\LR_{M}^{m}|f|(\om)
\leq&\lam_{0}^{-m}\LR_{A,\ph_{0}}^{m}|f|(\om)\\
\leq&\sum_{w\in S^{m}\,:\,w\cd\om_{0}\in W_{m+1}(A)}c_{\adr{Gibbsct1}}\mu_{0}[w]\left(\frac{1}{\mu_{0}[w]}\int_{[w]}f\,d\mu_{0}+[f]_{k}\theta^{m}\right)\\
\leq &c_{\adr{Gibbsct1}}\|f\|_{L^{1}(\mu_{0})}+c_{\adr{Gibbsct1}}[f]_{k}\theta^{m},
}
Moreover, for $\om,\up\in X$ with $\om_{0}\cdots\om_{k-1}=\up_{0}\cdots\up_{k-1}$
\ali
{
|\QR^{m}f(\om)-\QR^{m}f(\up)|\leq& \lam_{0}^{-m}\sum_{w\in S^{m}\,:\,w\cd \om_{0}\in W_{m+1}(M)}e^{S_{m}\ph(w\cd\om)}f(w\cd\om)-e^{S_{m}\ph(w\cd\up)}f(w\cd\up)|\\
\leq&\lam_{0}^{-m}\sum_{w\in S^{m}\,:\,w\cd \om_{0}\in W_{m+1}(M)}e^{S_{m}\ph(w\cd\om)}[f]_{k}\theta^{m} d_{\theta}(\om,\up)\\
+&\lam_{0}^{-m}\sum_{w\in S^{m}\,:\,w\cd \om_{0}\in W_{m+1}(M)}e^{S_{m}\ph(w\cd\up)}|f(w\cd\up)||e^{S_{m}\ph(w\cd\om)-S_{m}\ph(w\cd\up)}-1|\\
\leq&\theta^{m}\QR^{m}1(\om)[f]_{k}d_{\theta}(\om,\up)+(\QR^{m}|f|)(\up)e^{c_{\adr{Lamc}}d_{\theta}(\om,\up)}c_{\adr{Lamc}}d_{\theta}(\om,\up)\\
\leq&c_{\adr{Gibbsct1}}\theta^{m}[f]_{k}d_{\theta}(\om,\up)+(c_{\adr{Gibbsct1}}\|f\|_{L^{1}(\mu_{0})}+c_{\adr{Gibbsct1}}[f]_{k}\theta^{m})e^{c_{\adr{Lamc}}}c_{\adr{Lamc}}d_{\theta}(\om,\up).
}
Therefore the assertion of Claim \ref{cla:LYineqQ} is valid by putting $c_{\adr{Qct1}}=c_{\adr{Qct2}}=c_{\adr{Gibbsct1}}(1+2e^{c_{\adr{Lamc}}}c_{\adr{Lamc}})$.
\cla
{\label{cla:LYineqQ^k}
For any $\eta>0$ there exists $k_{0}\geq k$ such that $\|\QR^{k_{0}}f\|_{k}\leq c_{\adr{Qct1}}\|f\|_{L^{1}(\mu_{0})}+\|f\|_{k}(\theta(1+\eta))^{k_{0}}$ for any $f\in F^{k}_{b}(X)$.
}
This claim is satisfied by choosing $k_{0}$ in Claim \ref{cla:LYineqQ} so that $(c_{\adr{Qct2}})^{1/k_{0}}\leq \eta+1$.
\cla
{\label{cla:totalbdQ}
$\QR(\{f\in F^{k}_{b}(X)\,:\,\|f\|_{k}\leq 1\})$ is totally bounded in $\|\cd\|_{L^{1}(\mu_{0})}$.
}
It is sufficient to show that any sequence $f_{n}\in F^{k}_{b}(X)$ with $\|f_{n}\|_{k}\leq 1$ has a subsequence $\{f_{n(l)}\}_{l}$ so that $\QR f_{n(l)}$ converges in the sense of the norm $\|\cd\|_{L^{1}(\mu_{0})}$. 
By Claim \ref{cla:LYineqQ}, $\{\QR f_{n}(\om)\}$ is uniformly bounded in $\C$ for any $\om\in X$ and equicontinuous. 
Thus Ascoli Theorem and Lebesgue dominated convergence theorem tell us that such a sequence exists.
\cla
{\label{cla:bddL1Q}
There exists $c_{\adl{L1bddQ}}>0$ such that for any $f\in F^{k}_{b}(X)$, $\|\QR f\|_{L^{1}(\mu_{0})}\leq c_{\adr{L1bddQ}}\|f\|_{L^{1}(\mu_{0})}$.
}
Take the corresponding eigenfunction $h_{0}$ of the eigenvalue $\exp(P(\ph_{0}))$ of $\LR_{A,\ph_{0}}$. Then
\ali
{
\|\QR f\|_{L^{1}(\mu_{0})}
\leq\lam_{0}^{-1}\int_{X}\LR_{M}|f|\,d\mu_{0}
\leq&\lam_{0}^{-1}\int_{X}\LR_{A,\ph_{0}}|f|\,d\mu_{0}\\
\leq&\|h_{0}\|_{\infty}\lam_{0}^{-1}\int_{X}\LR_{A,\ph_{0}}|f|\,d\nu_{0}
\leq c_{\adr{L1bddQ}}\|f\|_{L^{1}(\mu_{0})}
}
with $c_{\adr{L1bddQ}}=\|h_{0}^{-1}\|_{\infty}\|h_{0}\|_{\infty}$. 

Hennion's theorem \cite[Theorem XIV.3]{HH} says that Claim \ref{cla:LYineqQ^k}, Claim \ref{cla:totalbdQ} and Claim \ref{cla:bddL1Q} imply that the essential spectral radius $r_{\text{ess}}(\QR)$ of $\QR\,:\,F^{k}_{b}(X)\to F^{k}_{b}(X)$ is less than or equal to $\theta(\eta+1)$. By arbitrary choosing $\eta>0$, we see $r_{\text{ess}}(\QR)\leq \theta$. From $\QR=\lam_{0}^{-1}\LR_{M}$, we obtain $r_{\text{ess}}(\LR_{M})\leq \lam_{0}\theta$. Hence the proof is complete.
\proe
Recall the quotient space $S/\!\!\leftrightarrow$ and the semi-order $\preceq$ on $S/\!\!\leftrightarrow$ defined in Section \ref{sec:Shiftsp}. 
For subset $T\subset S$, we write $\Si_{T}:=\bigcup_{a\in T}[a]$, and $\chi_{T}:=\chi_{\Si_{T}}$, the indicator of the set $\Si_{T}$. Denoted by $\rho(\LR)$ the resolvent set of $\LR\,:\,F^{k}_{b}(X)\to F^{k}_{b}(X)$.
\props
\label{prop:ex_efunc_geneRuelleop_noper_ex}
Assume that the conditions (A.1)-(A.3) are satisfied. Then
\ite
{
\item $r(\LR_{M})=\max_{T\in S/\leftrightarrow}r(\LR_{M(T)})$, where $M(T)$ means the submatrix of $M$ indexed by $T$;
\item $\bigcap_{T\in S/\leftrightarrow}\rho(\LR_{M(T)})\subset \rho(\LR_{M})$.
}
\prope
\pros
(1) Choose any $S_{1}\in S/\!\!\leftrightarrow$ and put $S_{2}=S\setminus S_{1}$. we let
$\LR_{ij}f=\chi_{S_{i}}\LR_{M}(\chi_{S_{j}}f)$
for $f\in F_{b}^{k}(X)$. 
By direct checking, we see $\|\LR_{ii}f\|_{\infty}\leq \|\LR_{M} f\|_{\infty}$ and $[\LR_{ii}f]_{k}\leq [\LR_{M}f]_{k}$. Therefore, we see
$\max_{i=1,2}r(\LR_{ii})\leq r(\LR_{M})$. Conversely, note that $\LR_{12}=\OR$ or $\LR_{21}=\OR$ by $S_{2}\not\preceq S_{1}$ or $S_{1}\not\preceq S_{2}$. We assume $\LR_{12}=\OR$. Then we obtain the expansion $\LR_{M}^{n}=\LR_{11}^{n}+\LR_{22}^{n}+\sum_{j=1}^{n}\LR_{22}^{j-1}\LR_{21}\LR_{11}^{n-j}$ for any $n\geq 1$. It is not hard to show that there exist constants $c_{\adl{cLii}}>0$ and $n_{0}\geq 1$ such that for any $n\geq n_{0}$, $\|\LR_{M}^{n}\|_{k}\leq c_{\adr{cLii}}\max_{i=1,2}r(\LR_{ii})^{n}$. Thus $r(\LR_{M})\leq \max_{i=1,2}r(\LR_{ii})$. The case $\LR_{21}=\OR$ is treated similarity. Thus we obtain $r(\LR_{M})= \max_{i=1,2}r(\LR_{ii})$.
\smallskip
\par
Choose any $0<\lam_{0}<r(\LR_{M})$. By the summability of $\ph$, there exists a decomposition $S/\!\!\leftrightarrow=\TR_{1}\cup \TR_{2}$ such that $\TR_{1}$ is finite and letting $T_{1}:=\bigcup_{T\in \TR_{2}}T$, $r(\LR_{M(T_{1})})$ is less than $\lam_{0}$. By using above argument repeatedly, we have $r(\LR_{M})=\max_{T\in S/\!\!\leftrightarrow}r(\LR_{M(T)})$ by the fact $r(\LR_{M})>r(\LR_{M(T_{1})})$. Hence the assertion is valid.
\smallskip
\\
(2) Take the notation $S_{1},S_{2},\LR_{ij}$ in (1). We display $\LR_{M}$ as the block operator matrix
\ali
{
\LR_{M}=\MatII{\LR_{11}&\LR_{12}\\\LR_{21}&\LR_{22}},
}
namely we decompose $\LR_{M}=\LR_{11}+\LR_{12}+\LR_{21}+\LR_{22}$ (see \cite{AJS,Tretter}).
Assume $\LR_{12}=\OR$. By Frobenius-Schur factorization \cite[Theorem 4.2]{AJS}, we have that for $\eta\in \rho(\LR_{11})$
\ali
{
\LR_{M}-\eta\IR=\MatII{\IR&\OR\\\LR_{21}(\LR_{11}-\eta\IR)^{-1}&\IR}\MatII{\LR_{11}-\eta\IR&\OR\\\OR&\LR_{22}-\eta\IR}.
}
Since the former $2\times 2$ block matrix is invertible and the inverse is bounded, $\eta\in \rho(\LR_{11})\cap \rho(\LR_{22})$ implies $\eta\in \rho(\LR_{M})$. In the case $\LR_{21}=\OR$, we obtain the same assertion by a similar argument. Choose any $\eta\in \bigcap_{T\in S/\leftrightarrow}\rho(\LR_{M(T)})$ and put $\lam_{0}=|\eta|$. We may assume $\eta\neq 0$ by $\ker \LR_{M(T)}\neq \{0\}$ for $T$ with $r(\LR_{M(T)})>0$. When we take $\TR_{1},\TR_{2}, T_{1}$ as the same notation in (1), we see $\eta\in \rho(\LR_{M(T_{1})})$ by $|\eta|>r(\LR_{M(T_{1})})$.
By the above argument repeatedly, we have $\rho(\LR_{M(T_{1})})\cap \bigcap_{T\in \TR_{1}}\rho(\LR_{M(T)})\subset \rho(\LR_{M})$. Since $\eta$ belongs to the left hand side, we obtain $\eta\in \rho(\LR_{M})$. Hence the proof is complete.
\proe
Finally we consider the case when $M$ is irreducible. 
We recall some notions in ergodic theory. A non-singular map $T$ on a sigma finite measure space $(X,\BR,\mu)$ is {\it conservative} if for every set $W\in \BR$ satisfying that $\{T^{-n}W\}$ are pairwise disjoint, $W = \emptyset $ or $X$ mod $\mu$. Such a set $W$ is called a {\it wandering set}.
The measure $\mu$ is called {\it exact} if $\mu(E)\mu(X\setminus E)=0$ for any $E\in \bigcap_{n=0}^{\infty}T^{-n}\BR$.
\props
\label{prop:ex_efunc_geneRuelleop_M_irre}
Assume that the matrix $M$ is irreducible with period $p$ and $\ph\,:\,X_{M}\to \R$ satisfies $[\ph]_{k+1}<\infty$ and $\|\LR_{M}1\|_{\infty}<\infty$. Assume also that there exist a triplet $(\lam,g,\nu)$ and an integer $0\leq i<p$ such that $\lam$ is equal to the spectral radius of $\LR_{M}\,:\,F^{k}_{b}(X_{M})\to F^{k}_{b}(X_{M})$, $g$ is a nonzero function in $\Lambda^{k}_{c}(X_{M}^{i})$ for some $c>0$ with $\LR_{M}^{p}g=\lam^{p} g$ on $X_{M}^{i}$, and $\nu$ is a nonzero Borel finite measure on $X^{i}_{M}$ with $\LR_{M}^{* p}\nu=\lam^{p} \nu$, where $X_{M}^{i}$ is defined in (\ref{eq:Sii=}). Then $\supp\,g=\supp\,\nu=X_{M}^{i}$ and putting $h=g/\nu(g)$, $\|\lam^{-np}\LR_{M}^{np}f-h\nu(f)\|_{L^{1}(X^{i}_{M},\nu)}\to 0$
as $n\to \infty$ for each $f\in L^{1}(X_{M}^{i},\nu)$.
\prope
\pros
The equalities of the support of $g$ and $\nu$ are proved by the proof of Proposition \ref{prop:supp_g} and Proposition \ref{prop:supp_nu}, respectively. Put $\mu=h\nu$. We define a bounded operator $\tLR f:=(\LR_{M}^{p}(h f))/(\lam^{p} h)$ for $f\in L^{1}(X_{M}^{i},\mu)$.
Consider the non-singular map $T:=\si^{p}$ on the measure space $(X_{M}^{i},\mu)$. Note that $\tLR$ becomes the transfer operator of $T$ on $L^{1}(X_{M}^{i},\mu)$ is given by
$\tLR f:=d\mu_{f}\circ T^{-1}/d\mu$ with $d\mu_{f}:=fd\mu$. Indeed, for $f_{1}\in L^{\infty}(X_{M}^{i})$ and $f_{2}\in L^{1}(X_{M}^{i})$
\ali
{
\int_{X_{M}^{i}}f_{1}\tLR f_{2}\,d\mu 
=\lam^{-p}\int_{X_{M}^{i}}\LR_{M}^{p}(f_{1}\circ T h f_{2})\,d\nu
=\int_{X_{M}^{i}}f_{1}\circ T f_{2}\,d\mu.
}
\par
First we will see that the measure $\mu$ is conservative. Choose any wandering set $W$. By the definition of $W$, we see $\sum_{n}\chi_{W}\circ T^{n}\leq 1$, where $\chi_{W}$ denotes the indicator of the set $W$. By Monotone convergence theorem,
\ali
{
\mu(1)\geq \int_{X_{M}^{i}}\sum_{n}\chi_{W}\circ T^{n}\,d\mu=\sum_{n=0}^{\infty}\int_{X_{M}^{i}}\chi_{W}\circ T^{n}\,d\mu=\sum_{n=0}^{\infty}\mu(W).
}
Since $\mu$ is finite, $\mu(W)$ must be zero. Thus $\mu$ is conservative.

Next we prove that the measure $\mu$ is exact. We define $\ti{S}=\{w=w_{0}\cdots w_{p-1}\in S_{0}^{p}\,:\,w_{0}\in S_{0,i} \text{ and }w\text{ is }M\text{-admissible}\}$ and a zero-one matrix $\ti{M}$ indexed as $\ti{S}$ by $\ti{M}(ww^\p)=1$ if $ww^\p$ is $M$-admissible and $\ti{M}(ww^\p)=0$ otherwise. Then we notice that $X_{M}^{i}$ is regarded as the topological Markov shift $X_{\ti{M}}$ with the state space $\ti{S}$, the transition matrix $\ti{M}$ and the shift $T$. Furthermore, $(X_{\ti{M}},T)$ is topological mixing, $\mu$ is a $T$-invariant Borel probability measure on $X_{\ti{M}}$, and the log Jacobian $\log(d\mu/d\mu\circ T)$ is equal to the potential $\ti{\ph}:=S_{p}\ph-\log h+\log h\circ T-\log\lam^p$. Observe that $\ti{\ph}\,:\,X_{\ti{M}}\to \R$ satisfies at least $[\ti{\ph}]_{k+1}<\infty$ with respect to $d_{\theta^{p}}$. Moreover, by recoding $X_{\ti{M}}$ below, the function $\ti{\ph}$ is deduced a function whose first variation is finite, i.e. $[\ti{\ph}]_{1}<\infty$. Indeed, let $S^{*}=\{w^{*}=w_{0}\cdots w_{k}\in \ti{S}^{k+1}\,:\,\ti{M}\text{-admissible}\}$. We give a zero-one matrix $M^{*}$ indexed as $S^{*}$ by $M^{*}(w^{*}v^{*})=1$ if and only if $w_{1}\cdots w_{k}=v_{0}\cdots v_{k-1}$ for $w^{*}=w_{0}\cdots w_{k}$ and $v^{*}=v_{0}\cdots v_{k}$. Then the TMS $X_{M^{*}}$ is conjugate to the TMS $X_{\ti{M}}$ by the map $\pi\,:\,X_{M^{*}}\to X_{\ti{M}}$, $\pi(w^{*}_{0}w^{*}_{1}\cdots)=w_{0}w_{1}\cdots$ for $w^{*}_{i}=w_{i}w_{i+1}\cdots w_{i+k}$. This map $\pi$ is bijective, bi-Lipschitz and satisfies that $\pi\circ T^{*}=T\circ \pi$ and $[\ti{\ph}\circ \pi]_{1}\leq [\ti{\ph}]_{k+1}<+\infty$, where $T^{*}$ is the shift-transformation on $X_{M^{*}}$. It is not hard to check that $\mu\circ \pi$ is a $T^{*}$-invariant Borel probability measure on $X_{M^{*}}$ and conservative. By virtue of \cite[Theorem 3.2]{ADU} (also \cite[Section 3.2]{Sar01_2} or \cite[Theorem 2.5]{Sar09}), it turns out that the measure $\mu\circ \pi$ is exact and so is $\mu$.

Finally we show the assertion of this proposition. By virtue of Lin's theorem (\cite{Lin} or \cite[Theorem 1.3.3]{A}), the exactness of $\mu$ means
\ali
{
\|\lam^{-np}\LR_{M}^{np}(hf_{0})-h\nu(hf_{0})\|_{L^{1}(X^{i}_{M},\nu)}
=&\|\tLR^{n}(f_{0}-\mu(f_{0}))\|_{L^{1}(X_{M}^{i},\mu)}\to 0
}
as $n\to \infty$ for each $f_{0}\in L^{1}(X_{M}^{i},\mu)$. Hence we obtain the assertion (2) by putting $f_{0}:=f/h\in L^{1}(X_{M}^{i},\mu)$ for any $f\in L^{1}(X_{M}^{i},\nu)$.
\proe
\section{Proofs of main results}\label{sec:proof}
\subsection{Proof of Theorem \ref{th:ex_efunc_geneRuelleop_M}(1)}\label{sec:proof_th:ex_efunc_geneRuelleop_M}
Assume that (A.1)-(A.3) and (B.1)-(B.3) are satisfied. In order to prove this theorem, we need to show the following lemmas.
We put $N:=\bigcup_{ij\,:\,A(ij)=1,M(ij)=0}[ij]$ and 
\alil
{
\ph(\e,\om):=
\case
{
\ph(\om)-(1/\e)\chi_{N}(\om),&\text{ if }\ph(\om)\geq \sup_{[\om_{0}]}\ph-1/\e\\
\sup_{[\om_{0}]}\ph-(1/\e)-(1/\e)\chi_{N}(\om),&\text{ if }\ph(\om)< \sup_{[\om_{0}]}\ph-1/\e\\
}\label{eq:ph-(1/e)chiN}
}
for $\e>0$ and $\om\in X$. Observe that $\ph(\om)\leq \ph(\e,\om)$ for $\om$ with $M(\om_{0}\om_{1})=1$.
Put $\psi:=\exp(\ph)(1-\chi_{N})$.
\lems\label{lem:seqpote}
Assume that (A.1)-(A.3) and (B.1)-(B.3) are satisfied. Then the potential (\ref{eq:ph-(1/e)chiN}) satisfies the condition (B.1)-(B.3).
\leme
\pros
We start with the validities of $[\phe]_{i}<\infty$ for $1\leq i\leq k$ and $[\phe]_{k+1}\leq [\ph]_{k+1}$. Let $\om,\up\in X$ with $\om_{0}\cdots \om_{l-1}=\up_{0}\cdots\up_{l-1}$ and $\om_{l}\neq \up_{l}$ for an integer $l\geq 1$. Note $\chi_{N}(\om)=\chi_{N}(\up)$ if $l\geq 2$. Consider the three cases:
\\
Case I: $\ph(\om),\ph(\up)\geq \sup_{[\om_{0}]}\ph-1/\e$. In this case, we have
\ali
{
|\ph(\e,\om)-\ph(\e,\up)|\leq& |\ph(\om)-\ph(\up)|+(1/\e)|\chi_{N}(\om)-\chi_{N}(\up)|\\
\leq&
\case
{
\sup_{[\om_{0}]}\ph-(\sup_{[\om_{0}]}\ph-1/\e)+1/\e,& \text{ if }1\leq l\leq k\\
|\ph(\om)-\ph(\up)|,& \text{ if }l> k\\
}\\
\leq& 
\case
{
(2/(\theta^{k}\e))d_{\theta}(\om,\up),& \text{ if }1\leq l\leq k\\
[\ph]_{k+1}d_{\theta}(\om,\up),& \text{ if }l> k.\\
}
}
Case II: $\ph(\up)< \sup_{[\om_{0}]}\ph-1/\e\leq\ph(\om)$ or $\ph(\om)< \sup_{[\om_{0}]}\ph-1/\e\leq\ph(\up)$. We obtain
\ali
{
|\ph(\e,\om)-\ph(\e,\up)|\leq& \max(\ph(\om),\ph(\up))-(\textstyle\sup_{[\om_{0}]}\ph-1/\e)+(1/\e)|\chi_{N}(\om)-\chi_{N}(\up)|\\
\leq& 
\case
{
2/\e,& \text{ if }1\leq l\leq k\\
\max(\ph(\om),\ph(\up))-\min(\ph(\om),\ph(\up)),& \text{ if }l> k.\\
}
}
Therefore the same assertion as in Case I is satisfied.
\\
\smallskip
Case III: $\ph(\om),\ph(\up)< \sup_{[\om_{0}]}\ph-1/\e$. We see $|\ph(\e,\om)-\ph(\e,\up)|\leq 1/\e$ if $l=1$ and equals $0$ if $l\geq 2$.
Thus the condition (B.1) is fulfilled. The condition (B.2) is valid from $\sup_{\e>0}\sup_{\om\in [s]}\ph(\e,\om)\leq \sup_{\om\in [s]}\ph$ for any $s\in S$ and the summability of $\ph$.
The condition (B.3) follows from for each $a\in S$ and $\om\in [a]$,
\ali
{
|e^{\ph(\e,\om)}-\psi(\om)|\leq&
\case
{
e^{\sup_{[\om_{0}]}}e^{-1/\e}, &\text{ if }M(\om_{0}\om_{1})=0\\
0, &\text{ if }M(\om_{0}\om_{1})=1 \text{ and }\ph(\om)\geq \sup_{[\om_{0}]}\ph-1/\e\\
e^{\sup_{[\om_{0}]}\ph-1/\e}+e^{\ph(\om)}, &\text{ if }M(\om_{0}\om_{1})=1 \text{ and }\ph(\om)< \sup_{[\om_{0}]}\ph-1/\e\\
}\\
\leq&2e^{\sup_{[a]}}e^{-1/\e}\to 0.
}
Hence (B.1)-(B.3) are fulfilled.
\proe
\lems\label{lem:ext_efev}
Assume that (A.1)-(A.3) and (B.1)-(B.3) are satisfied. Let $\lam:=r(\LR_{M})$. Then there exist a nonnegative function $g\in \Lambda^{k}_{c_{\adr{Lamc}}}$ with $\|g\|_{\infty}=1$ and the Borel probability measure $\nu$ on $X$ such that $\LR_{M}g=\lam g$ and $\LR_{M}^{*}\nu=\lam \nu$.
\leme
\pros
First we assume the condition $(A.1)^{\p}$, i.e. $A$ is finitely irreducible. By this assumption and by Lemma \ref{lem:seqpote}, Theorem \ref{th:exGibbs} implies that there exists a triplet $(\lam(\e),g(\e,\cd),\nu(\e,\cd))\in \R\times \Lambda^{k}_{c_{\adr{udphe2}}}\times C_{b}(X)^{*}$ such that $\lam(\e)=\exp(P(\phe))$
\alil
{
&\LR_{A,\phe}h(\e,\cd)=\lam(\e)h(\e,\cd),\ \  \LR_{A,\phe}^{*}\nu(\e,\cd)=\lam(\e)\nu(\e,\cd),\ \ \|g(\e,\cd)\|_{\infty}=\nu(\e,1)=1,\label{eq:pTSC}
}
where $g(\e,\cd)$ is defined by $h(\e,\cd)/\|h(\e,\cd)\|_{\infty}$.
By noting that $\e \mapsto P(\ph(\e,\cd))$ is increasing, the limit $\lam_{1}:=\lim_{\e\to 0}\exp(P(\ph(\e,\cd))$ exists.
In addition to the fact $c_{\adr{udphe}}\leq [\ph]_{k+1}$, Proposition \ref{prop:conv_efunc_LB} implies that there exists a nonnegative function $g\in \Lambda^{k}_{c_{\adr{Lamc}}}$ such that $\LR_{M}g=\lam_{1} g$ and $\|g\|_{\infty}=1$. Therefore $\lam_{1}\leq \lam:=r(\LR_{M})$. Moreover by
\ali
{
\|\LR_{A,\phe}^{n}1\|_{k}\geq\|\LR_{A,\phe}^{n}1\|_{\infty}\geq \|\LR_{M,\phe}^{n}1\|_{\infty}\geq\|\LR_{M,\ph}^{n}1\|_{\infty}=\|\LR_{M}^{n}\|_{\infty},
}
we notice $\lam_{1}\geq \lam$ and thus $\lam=\lam_{1}$. 
On the other hand, Proposition \ref{prop:conv_evec_LB} yields the existence of a Borel probability measure $\nu$ on $X$ such that $\LR_{M}^{*}\nu=\lam\nu$.
\smallskip
\par
Next we consider the general case (A.1). We take a finitely irreducible matrix $\hat{A}$ indexed by $S$ so that $A(ij)=1$ implies $\hat{A}(ij)=1$. For example, we let $\hat{A}(ij)\equiv 1$. Fix $\e>0$. By Proposition \ref{prop:ext_pote}, there exists $\hat{\ph}(\e,\cd)\,:\,X_{\hat{A}}\to \R$ such that $\hat{\ph}(\e,\cd)=\ph(\e,\cd)$ on $X_{A}$,  $[\hat{\ph}(\e,\cd)]_{k+1}=[\ph(\e,\cd)]_{k+1}$ and $\hat{\ph}(\e,\cd)$ is summable. By replacing A and M with $\hat{A}$ and $A$, respectively and repeating the argument above, we obtain a triplet $(\hat{\lam}(\e),\hat{g}(\e,\cd),\hat{\nu}(\e,\cd))$ satisfying that $\hat{\lam}(\e)=r(\LR_{A,\hat{\ph}(\e,\cd)})$, $\hat{g}(\e,\cd)\in \Lambda^{k}_{c_{\adr{Lamc}}}$, $\hat{\nu}(\e,\cd)$ is a Borel probability measure on $X_{\hat{A}}$ and
\ali
{
&\LR_{A,\hat{\ph}(\e,\cd)}\hat{g}(\e,\cd)=\hat{\lam}(\e)\hat{g}(\e,\cd),\ \  \LR_{A,\hat{\ph}(\e,\cd)}^{*}\hat{\nu}(\e,\cd)=\hat{\lam}(\e)\hat{\nu}(\e,\cd),\ \ \|\hat{g}(\e,\cd)\|_{\infty}=\hat{\nu}(\e,1)=1.
}
According to Proposition \ref{prop:supp_nu} and Proposition \ref{prop:supp_g}, we see that $\supp\, \hat{g}(\e,\cd)=X_{\hat{A}}$ and $\supp\, \hat{\nu}(\e,\cd))=X_{A}$. Then we define $g(\e,\cd):=\hat{g}(\e,\cd)\chi_{X_{A}}/\|\hat{g}(\e,\cd)\chi_{X_{A}}\|_{\infty}$ and $\nu(\e,\cd):=\nu(\e,\cd|X_{A})$. By considering the restricted operator $\LR_{A,\ph(\e,\cd)}=\LR_{A,\hat{\ph}(\e,\cd)}|_{F_{b}^{k}(X_{A})}$ on $F_{b}^{k}(X_{A})$, we get the same equations as (\ref{eq:pTSC}). The equation $\hat{\lam}(\e)=\lam(\e)$ is guaranteed from $\hat{\lam}(\e)\geq \hat{\lam}(\e)$ by the definition and from $\hat{\lam}(\e)$ is an eigenvalue of $\LR_{A,\ph(\e,\cd)}$. By a similar argument above again, we obtain the assertion. Hence the proof is complete.
\proe
\lems\label{lem:lam=ePph}
Assume that (A.1)-(A.3) and (B.1)-(B.3) are satisfied. Then $\lam=\exp(P(\ph))$.
\leme
\pros
In view of Proposition \ref{prop:ex_efunc_geneRuelleop_noper_ex}(1), there exists $T\in S/\!\!\leftrightarrow$ such that $r(\LR_{M(T)})=\lam$. By using Proposition \ref{prop:conv_efunc_LB} and proposition \ref{prop:supp_g} replacing $M$ by $M(T)$, we get the corresponding eigenfunction $g_{0}\in \Lambda^{k}_{c_{\adr{Lamc}}}$ of $\LR_{M(T)}$ with $\supp\,g_{0}=\bigcup_{a\in T}[a]$. Take any $\om\in X_{M(T)}$.
We see
\ali
{
\lam^{n}g_{0}(\om)=&\LR_{M(T)}^{n}g_{0}(\om)\leq \LR_{M}^{n}g_{0}(\om)\leq \sum_{w\in W_{n}(M)}\exp(\sup_{\up\in [w]\cap X_{M}}S_{n}\ph(\up)).
}
Therefore $\log \lam\leq P(\ph|_{X_{M}})$. On the other hand, in addition to the fact $\ph(\om)\leq \ph(\e,\om)$ for $\om$ with $M(\om_{0}\om_{1})=1$, we have that for $w\in W_{n}(M)$
\ali
{
\exp(\sup_{\om\in [w]\cap X_{M}}S_{n}\ph(\om))\leq \exp(\sup_{\om\in [w]\cap X_{M}}S_{n}\ph(\e,\om))
\leq \exp(\sup_{\om\in [w]}S_{n}\ph(\e,\om)).
}
Thus $P(\ph|_{X_{M}})\leq P(\ph(\e,\cd))$ and $P(\ph|_{X_{M}})\leq \log \lam$. Hence the assertion is fulfilled.
\proe
\pros[Proof of Theorem \ref{th:ex_efunc_geneRuelleop_M}(1)] It follows from Lemma \ref{lem:ext_efev} and Lemma \ref{lem:lam=ePph}.
\proe
The proof of Theorem \ref{th:ex_efunc_geneRuelleop_M}(2) will be carried over after the proof of Theorem \ref{th:ex_efunc_geneRuelleop_M_irre}.
\subsection{Proof of Theorem \ref{th:ex_efunc_geneRuelleop_M_irre}}\label{sec:proof_th:ex_efunc_geneRuelleop_M_irre}
First we prove this theorem under the assumption $(A.1)^{\p}$.
\lems\label{lem:ress<=thetalam}
Assume that $(A.1)^\p$, (A.2) and (A.3) are satisfied. Then $r_{\text{ess}}(\LR_{M})\leq \theta r$. 
\leme
\pros
By virtue of Proposition \ref{prop:conv_quasicomp}, we get the inequality $r_{\text{ess}}(\LR_{M})\leq \exp(P(\phe))\theta$ for any $\e>0$. Letting $\e\to 0$, we obtain $r_{\text{ess}}(\LR_{M})\leq \lam\theta$.
\proe
\lems\label{lem:ex_efunc_geneRuelleop_M_irre_fi}
Assume that $(A.1)^\p$, (A.2) and (A.3) are satisfied. Then the assertion of Theorem \ref{th:ex_efunc_geneRuelleop_M_irre} holds.
\leme
To show this, we start with the general spectral form of $\LR_{M}\,:\,F^{k}_{b}(X)\to F^{k}_{b}(X)$. By virtue of Lemma \ref{lem:ress<=thetalam}, we have the form
\alil
{
\textstyle\LR_{M}=\sum_{j=0}^{q-1}(\lam_{j}\PR_{j}+\DR_{j})+\RR\label{eq:LMph=...}
}
satisfying the following (1)-(5):
\ite
{
\item[(1)] each $\lam_{j}\in \C$ are eigenvalues of $\LR_{M}$ with finite multiplicity and with $|\lam_{j}|=\lam$;
\item[(2)] each $\PR_{j}$ is the projection onto the generalized eigenspace associated to $\lam_{j}$;
\item[(3)] each $\DR_{j}$ is nilpotent, i.e. $\DR_{j}^{n_{j}-1}\neq \OR$ and $\DR_{j}^{n_{j}}=\OR$ for some $n_{j}\geq 1$; 
\item[(4)] $\PR_{i}\LR_{M}=\LR_{M}\PR_{i}$, $\PR_{i}^{2}=\PR_{i}$, $\PR_{i}\DR_{i}=\DR_{i}\PR_{i}=\DR_{i}$ for $i$ and $\PR_{i}\PR_{j}=\DR_{i}\DR_{j}=\PR_{i}\RR=\RR\PR_{i}=\RR\DR_{i}=\DR_{i}\RR=\OR$ for each $i\neq j$;
\item[(5)] the spectral radius of $\RR$ is less than $\lam$.
}
\par
Take the eigenfunction $g$ and the eigenvector $\nu$ of the eigenvalue $\lam$ given in Theorem \ref{th:ex_efunc_geneRuelleop_M}(1). Since $M$ is irreducible, the equalities $\supp\, g=\bigcup_{s\in S_{0}}[s]$ and $\supp\,\nu=X_{M}$ hold by Proposition \ref{prop:supp_g} and Proposition \ref{prop:supp_nu}.
We need the following ten claims:
\ncla
{\label{cla:semiform1}
$\LR_{M}^{n}f=\sum_{i=0}^{q-1}\lam_{i}^{n}\PR_{i}f+\RR^{n} f \text{ on }X_{M}$ for any $f\in F^{k}_{b}(X)$ and $n\geq 1$.
}
Choose any $f\in F^{k}_{b}(X)$ and $0\leq i\leq q-1$. We will show $\DR_{i}^{s}f=0$ on $X_{M}$ for each $s=n_{i}-1,n_{i}-2,\dots, 1$ inductively. We may assume $n_{i}\geq 2$. Put $\xi=\DR_{i}^{s-2}\PR_{i}f$ with $s=n_{i}-1$. Then $\DR_{i}\xi=\DR_{i}^{s}f$ and $\DR_{i}^{l}\xi=0$ on $X_{M}$ for $l\geq 2$. Note also $\PR_{i}\xi=\xi$. The equation (\ref{eq:LMph=...}) implies
\alil
{
\LR_{M}^{n}\xi=&\sum_{j=0}^{q-1}\Big(\lam_{j}^{n}\PR_{j}\xi+\sum_{s=1}^{n_{i}-1}\binom{n}{s}\lam_{j}^{n-s}\DR_{j}^{s}\xi\Big)+\RR^{n}\xi=\lam_{i}^{n}\xi+n\lam_{i}^{n-1}\DR_{i}\xi\label{eq:L^nxi=...}
}
on $X_{M}$. Therefore
\alil
{
\|\lam^{-n}\LR_{M}^{n}\xi\|_{L^{1}(\nu)}\geq n\lam^{-1}\|\DR_{i}\xi\|_{L^{1}(\nu)}-\|\xi\|_{L^{1}(\nu)}.\label{eq:Ln>=n...}
}
The left hand side is bounded by $\nu(|\xi|)$.
Thus $\nu(|\DR_{i}\xi|)$ must be zero. Since $\nu$ is positive with full support on $X_{M}$ and $\DR_{i}\xi$ is continuous, we get $\DR_{i}\xi=\DR_{i}^{s}f=0$ on $X_{M}$.

Assume $\DR_{i}^{j}f=0$ on $X_{M}$ for $j=n_{i}-1,n_{i}-2,\dots, s$ with $s\geq 2$ and put $\xi=\DR_{i}^{j-1}\PR_{i}f$ with $j=s-1$. By repeatedly considering (\ref{eq:L^nxi=...}) and (\ref{eq:Ln>=n...}), we get the equation $\DR_{i}\xi=\DR_{i}^{j}f=0$ on $X_{M}$. Hence we obtain the claim.
\cla
{\label{cla:bddonXM2}
$c_{\adl{bdLM}}:=\sup_{n\geq 1}\|\lam^{-n}\LR_{M}^{n}1\|_{\infty}<+\infty$.
}
Recall $\LR_{M}^{n}1\in \Lambda_{c_{\adr{Lamc}}}^{k}$ for $n\geq 1$ by Proposition \ref{prop:prop_Lamc0}(3). For any $\om\in X$ with $\om_{0}\in S_{0}$, we have that by taking $\up\in X_{M}$ with $\om_{0}=\up_{0}\in S_{0}$, 
\ali
{
\lam^{-n}\LR_{M}^{n}1(\om)=\lam^{-n}\LR_{M}^{k}(\LR_{M}^{n-k}1)(\om) 
\leq&e^{c_{\adr{Lamc}}\theta^{k}}\lam^{-n}\sum_{w\in S_{0}^{k}\,:\,w\cd\om_{0}\in W_{k+1}(M)}e^{S_{k}\ph(w\cd\om)}\LR_{M}^{n-k}1(w\cd\up)\\
\leq&e^{c_{\adr{Lamc}}\theta^{k}}\lam^{-k}\|\LR_{M}^{k}1\|_{\infty}(\sum_{i=0}^{q-1}\|\PR_{i}1\|_{\infty}+\sup_{j}\|\lam^{-j}\RR^{j}1\|_{k})
}
for any $n>k$ by using Claim \ref{cla:semiform1}. Since the last expression is finite, the claim is satisfied.
\cla
{\label{cla:semisimple}
The operator $\LR_{M}$ has the form $\LR_{M}=\sum_{i=0}^{q-1}\lam_{i}\PR_{i}+\RR$.
}
Suppose $n_{i}\geq 2$, i.e. $\DR_{i}^{n_{i}-1}\neq \OR$ and $\DR_{i}^{n_{i}}=\OR$. Let $f\in F^{k}_{b}(X)$ and $\xi=\DR_{i}^{s-1}\PR_{i}f$ with $s=n_{i}-1$. Observe $\DR_{i}\xi=\DR_{i}^{s}f$. It follows from (\ref{eq:L^nxi=...}) and the claim \ref{cla:bddonXM2} that $n\lam^{-1}\|\DR_{i}\xi\|_{\infty}\leq c_{\adr{bdLM}}+\|\xi\|_{\infty}<+\infty$ for any $n$. Therefore $\|\DR_{i}\xi\|_{\infty}=\|\DR_{i}^{s}f\|_{\infty}=0$. By a similar argument for each $s=n_{i}-2,\dots, 1$ repeatedly, we obtain the assertion.
\cla
{\label{cla:basis_nonfunc}
If $f\in F^{k}_{b}(\Si_{j})$ satisfies $\LR_{M}^{p}f=\lam^{p}f$ on $\Si_{j}$ for some $j$, then there exist nonnegative functions $f_{i}\in F^{k}_{b}(\Si_{j})$ $(i=0,1,2,3)$ such that $f=f_{0}-f_{1}+\sqrt{-1}f_{2}-\sqrt{-1}f_{3}$ on $\Si_{j}$ and $\LR_{M}^{p}f_{i}=\lam^{p}f_{i}$ on $\Si_{j}$ for $i=0,1,2,3$.
}
Since $\LR_{M}^{p}\Re f=\lam^{p}\Re f$ and $\LR_{M}^{p}\Im f=\lam^{p}\Im f$ on $\Si_{j}$, we may assume $f=\Re f$ and $f=0$ on $X\setminus \Si_{j}$.
By Claim \ref{cla:bddonXM2} and the basic inequality (\ref{eq:basine}), $c_{\adl{bffL^n}}:=\sup_{n\geq 1}\|\lam^{-n}\LR_{M}^{n}\|_{k}$ is finite. We decompose into $f=f_{+}-f_{-}$ with $f_{+}\geq 0$ and $f_{-}\geq 0$. Since $|f_{\pm}(\om)-f_{\pm}(\up)|\leq |f(\om)-f(\up)|$ holds, we get $f_{+},f_{-}\in F^{k}_{b}(X)$. Note the form
$f=(1/n)\sum_{i=0}^{n-1}\lam^{-ip}\LR_{M}^{ip}f
=f^{(n)}_{+}-f^{(n)}_{-}$
by putting $f^{(n)}_{\pm}=(1/n)\sum_{ki0}^{n-1}\lam^{-ip}\LR_{M}^{ip}f_{\pm}$. In addition to the fact $\|f^{(n)}_{\pm}\|_{k}\leq c_{\adr{bffL^n}}\|f\|_{k}$ for any $n\geq 1$, Ascoli Theorem implies that there exist a subsequence $(n_{i})$ and functions $f_{\pm}^{(\infty)}\,:\,X\to \R$ such that $f_{\pm}^{(n_{i})}(\om)\to f_{\pm}^{(\infty)}(\om)$ as $i\to \infty$ for each $\om\in X$. By $\sup_{n\geq 1}[f_{\pm}^{(n)}]_{k}<\infty$, we see $f_{+}^{(\infty)}, f_{-}^{(\infty)}\in F^{k}_{b}(X)$. Take any potential $\ti{\ph}\in F^{1}(X,\R)$ satisfying that $\ti{\lam}=\exp(P(\ti{\ph}))<\infty$ and $\LR_{M}\xi\leq \LR_{A,\ti{\ph}}\xi$ for any $\xi\geq 0$ (e.g. (\ref{eq:ph-(1/e)chiN}) for the existence). Let $\ti{\mu}$ be the corresponding positive eigenvector of $\ti{\lam}$ of $\LR_{A,\ti{\ph}}^{*}$ with full measure on $X$ and $\ti{\mu}(1)=1$. We have
$\lam^{-p}\LR_{M}^{p}f_{\pm}^{(n)}-f_{\pm}^{(n)}=(1/n)(\lam^{-np-p}\LR_{M}^{np+p}f_{\pm}-f_{\pm})\to 0$ in $\|\cd\|_{k}$, and
\ali
{
\|\LR_{M}^{p} f_{\pm}^{(n_{i})}-\LR_{M}^{p}f_{\pm}^{(\infty)}\|_{L^{1}(\ti{\mu})}
\leq \int_{X}\LR_{M}^{p}|f_{\pm}^{(n_{i})}-f_{\pm}^{(\infty)}|\,d\ti{\mu}
\leq&\int_{X}\LR_{A,\ti{\ph}}^{p}|f_{\pm}^{(n_{i})}-f_{\pm}^{(\infty)}|\,d\ti{\mu}\\
=&\ti{\lam}^{p}\|f_{\pm}^{(n_{i})}-f_{\pm}^{(\infty)}\|_{L^{1}(\ti{\mu})}\to 0
}
as $i\to \infty$ using Lebesgue dominated convergence theorem. Consequently 
we get $f_{\pm}^{(\infty)}=\lam^{-p}\LR_{M}^{p}f_{\pm}^{(\infty)}$ $\ti{\mu}$-a.e. By continuity, we obtain $\LR_{M}^{p}f_{\pm}^{(\infty)}=\lam^{p} f_{\pm}^{(\infty)}$ on $X$. 
By noting the form $f=f_{+}^{(\infty)}-f_{-}^{(\infty)}$, the assertion of the claim is yielded.
\cla
{\label{cla:inLam}
If $f\in F^{k}_{b}(\Si_{i})$ satisfies $\LR_{M}^{p}f= \lam^{p} f$ and $f\geq  0$ on $\Si_{i}$, then $f\in \Lambda^{k}_{c_{\adr{Lamc}}}(\Si_{i})$.
}
For any $\om,\up \in \Si_{i}$ with $d_{\theta}(\om,\up)\leq \theta^{k}$ and $n\geq 0$ with $np>k$, we have
\ali
{
f(\om)=\lam^{-np}\LR_{M}^{np}f(\om)
\leq &\lam^{-np}\sum_{w\in S_{0}^{n}\,:\,w\cd\om_{0}\in W_{np+1}(M)}e^{S_{np}\ph(w\cd\up)+c_{\adr{Lamc}}d_{\theta}(\om,\up)}(f(w\cd\up)+[f]_{k}\theta^{np+k})\\
\leq& e^{c_{\adr{Lamc}}d_{\theta}(\om,\up)}f(\up)+c_{\adr{bffL^n}}[f]_{k}\theta^{np}.
}
Letting $n\to \infty$, we see $f(\om)\leq e^{c_{\adr{Lamc}}d_{\theta}(\om,\up)}f(\up)$ and thus $f\in \Lambda^{k}_{c_{\adr{Lamc}}}(\Si_{i})$.
\cla
{\label{cla:simple_|f|}
If $f\in F^{k}_{b}(X)$ and $\eta\in \C$ satisfy $\LR_{M}f=\eta f$ and $|\eta|=\lam$, then for each $i=0,1,\dots, p-1$ there exists a constant $c\geq 0$ such that 
$|f|=c h$ on $X_{M}^{i}$.
}
By $\lam|f|=|\eta f|=|\LR_{M}f|\leq \LR_{M}|f|$ and $\nu(\LR_{M}|f|-\lam|f|)=0$, we have $\LR_{M}|f|=\lam|f|$ $\nu$-a.e. From $\supp\,\nu=X_{M}$, we see $\LR_{M}|f|=\lam|f|$ on $X_{M}$ and therefore $\LR_{M}^{p}|f|=\lam^{p}|f|$ on $X_{M}^{i}$. Put $\ti{\nu}:=\nu(\cd|X_{M}^{i})$ and $\ti{h}:=h\chi_{\Si_{i}}/\ti{\nu}(h)$. By virtue of Proposition \ref{prop:ex_efunc_geneRuelleop_M_irre}, we see
\ali
{
\||f|-\ti{h}\ti{\nu}(|f|)\|_{L^{1}(\ti{\nu})}=
\|\lam^{-np}\LR_{M}^{np}|f|-\ti{h}\ti{\nu}(|f|)\|_{L^{1}(\ti{\nu})}\to 0
}
as $n\to \infty$. 
By continuity of $|f|$ and $\hat{h}$, it turns out that $|f|=c h$ on $X_{M}^{i}$ with $c=\hat{\nu}(|f|)/\hat{\nu}(h)$. Hence the claim is fulfilled.
\smallskip
\par
Recall $\lam_{i}=\lam\kappa^{i}$, $h_{i}=\sum_{j=0}^{p-1}\kappa^{-ji}h\chi_{\Si_{j}}$ and $\nu_{i}=\sum_{j=0}^{p-1}\kappa^{ji}\nu|_{\Si_{i}}$ in Theorem \ref{th:ex_efunc_geneRuelleop_M_irre} with $\kappa=\exp(2\pi\sqrt{-1}/p)$. In addition to the fact $\LR_{M}(\chi_{\Si_{i}}f)=\chi_{\Si_{i+1}}\LR_{M}f$ for $f\in C_{b}(X)$, it follows that the equations $\LR_{M}h_{i}=\lam_{i}h_{i}$ and $\LR_{M}^{*}\nu_{i}=\lam_{i}\nu_{i}$.
\cla
{\label{cla:lamk_simple}
For $0\leq l<p$, the eigenvalue $\lam_{l}$ is simple, namely if $f\in F^{k}_{b}(X)$ satisfies $\LR_{M}f=\lam_{l} f$ then $f=c h_{l}$ for some constant $c\in \C$.
}
Choose any $0\leq j<p$. By noting the equation $\LR_{M}^{p}f=\lam^{p} f$ on $\Si_{j}$, it follows from Claim \ref{cla:basis_nonfunc} and Claim \ref{cla:inLam} that there exist $f_{i}\in \Lambda^{k}_{c_{\adr{Lamc}}}(\Si_{j})$ $(i=0,1,2,3)$ such that $f=f_{0}-f_{1}+\sqrt{-1}(f_{2}-f_{3})$ on $\Si_{j}$, $\LR_{M}^{p}f_{i}=\lam^{p} f_{i}$ on $\Si_{j}$. Claim \ref{cla:simple_|f|} tells us that equation $f_{i}=c_{\adr{csi}}(ij)h$ on $X_{M}^{j}$ for some constant $c_{\adl{csi}}(ij)\geq 0$. Moreover, this equation is extended on $\Si_{j}$. Indeed, for any $\om,\up\in \Si_{j}$ with $\up\in X_{M}$ and $\om_{0}=\up_{0}$ and for any $w\in S_{0}^{n}$ with $w\cd\om_{0}\in W_{n+1}(M)$, we notice $f_{i}(w\cd\om)\leq e^{c_{\adr{Lamc}}\theta^{n}}f_{i}(w\cd\up)=e^{c_{\adr{Lamc}}\theta^{n}}c_{\adr{csi}}(ij)h(w\cd\up)\leq e^{2c_{\adr{Lamc}}\theta^{n}}c_{\adr{csi}}(ij)h(w\cd\om)$. This observation implies for any $n\geq k$
\ali
{
f_{i}(\om)=\lam^{-n}\LR_{M}^{n}f_{i}(\om)\leq e^{2c_{\adr{Lamc}}\theta^{n}}\lam^{-n}\LR_{M}^{n}(c_{\adr{csi}}(ij)h(\om))=e^{2c_{\adr{Lamc}}\theta^{n}}c_{\adr{csi}}(ij)h(\om).
}
Letting as $n\to \infty$, we get $f_{i}(\om)\leq c_{\adr{csi}}(ij)h(\om)$. Similarity, the converse $c_{\adr{csi}}(ij)h(\om)\leq f_{i}(\om)$ holds and thus $f_{i}=c_{\adr{csi}}(ij)h$ on $\Si_{j}$. Consequently, we obtain $f=c_{\adr{csi2}}(j)h$ on $\Si_{j}$ by putting $c_{\adl{csi2}}(j)=c_{\adr{csi}}(0j)-c_{\adr{csi}}(1j)+\sqrt{-1}(c_{\adr{csi}}(2j)-c_{\adr{csi}}(3j))$. Finally we prove the equation $f=c_{\adr{csi2}}(0)h_{l}$ on $X$. We notice the equation
$\lam_{l}c_{\adr{csi2}}(j)h\chi_{\Si_{j}}=\lam_{l}f\chi_{\Si_{j}}=\LR_{M}(\chi_{\Si_{j-1}}f)=\chi_{\Si_{j}}c_{\adr{csi2}}(j-1)\lam h$.
Therefore $c_{\adr{csi2}}(j)=\lam(\lam_{l})^{-1}c_{j-1}=\kappa^{-l}c_{j-1}=\kappa^{-jl}c_{\adr{csi2}}(0)$. Thus $f=c_{\adr{csi2}}(0)h_{l}$ is valid by the definition of $h_{l}$.
\cla
{\label{cla:|f|onXM}
If $\LR_{M}f=\eta f$ for some $f\in F^{k}_{b}(X)$ and $|\eta|=\lam$, then
$\|f\|_{\infty}\leq c_{\adr{bffL^n}}\sup_{\om\in X_{M}}|f(\om)|$.
}
Let $\om\in \supp f\subset \bigcup_{i}\Si_{i}$ and
take $\up\in X_{M}$ so that $\om_{0}=\up_{0}$. We have
\ali
{
|f|(\om)\leq& \lam^{-n}\LR_{M}^{n}|f|(\om)\\
\leq&\lam^{-n}\sum_{w\,:\,w\cd\om_{0}\in W_{n+1}(M)}e^{S_{n}\ph(w\cd\om)}(|f|(w\cd\up)+[f]_{k}\theta^{n})
\leq c_{\adr{bffL^n}}\sup_{X_{M}}|f|+c_{\adr{bffL^n}}[f]_{k}\theta^{n}.
}
for any $n>k$. Letting $n\to \infty$, we obtain the assertion.
\cla
{\label{cla:p=q}
If $f\in F^{k}_{b}(X)$ and $\eta\in \C$ satisfy $\LR_{M}f=\eta f$, $f\neq 0$ and $|\eta|=\lam$, then $\eta^{p}=\lam^{p}$. Namely, the number $q$ in Claim \ref{cla:semisimple} equals $p$.
}
By virtue of Claim \ref{cla:|f|onXM}, we notice $f\neq 0$ on $X_{M}^{i}$ for some $i$.
Put $\ti{\nu}=\nu(\cd|X_{M}^{i})$. 
Proposition \ref{prop:ex_efunc_geneRuelleop_M_irre} implies that  $\eta^{np}\lam^{-np}f=\lam^{-np}\LR_{M}^{pn}f$ converges in $L^{1}(\ti{\nu})$ as $n\to \infty$. By $\|f\|_{L^{1}(\ti{\nu})}>0$, the number $\eta^{p}\lam^{-p}$ must be $1$. 
\cla
{\label{cla:Pi=hinui}
$P_{i}=h_{i}\otimes \nu_{i}$ for each $i$.
}
Since $\lam_{i}$ is simple, for any $f\in F^{k}_{b}(X)$ there exists a unique number $\tau(f)\in \C$ such that $\PR_{i} f=\tau(f)h_{i}$. It is no hard to check that $\tau$ is a linear functional and bounded by $\|\PR_{i}\|_{k}/\|h_{i}\|_{k}$. Since $\tau(\LR_{M}f)h_{i}=\PR_{i}(\LR_{M}f)=\lam_{i}\PR_{i}f=\lam_{i}\tau(f)h_{i}$ holds, $\tau$ satisfies $\LR_{M}^{*}\tau=\lam_{i}\tau$.
Note the equation
$\PR_{i}f=\lim_{n\to \infty}(1/n)\sum_{j=0}^{n-1}\lam_{i}^{-j}\LR_{M}^{j}f$ in $F^{k}_{b}(X)$ (e.g. \cite[Corollary III.4]{HH}). Therefore we have
$\nu_{i}(f)=\nu_{i}((1/n)\sum_{j=0}^{n-1}\lam_{i}^{-j}\LR_{M}^{j}f)\to \nu_{i}(\PR_{i}f)=\tau(f)\nu_{i}(h_{i})=\tau(f)$
as $n\to \infty$. Thus $\nu_{i}=\tau$ holds.
\smallskip
\par
Hence Lemma \ref{lem:ex_efunc_geneRuelleop_M_irre_fi} follows from Claim \ref{cla:semisimple}, Claim \ref{cla:lamk_simple} and Claim \ref{cla:Pi=hinui}. 
\qed
\pros[Proof of Theorem \ref{th:ex_efunc_geneRuelleop_M_irre}]
Assume that (A.1)-(A.3) are satisfied. Take any topological Markov shift $X_{\hat{A}}$ whose transition matrix is $S\times S$ finitely primitive and satisfies $A(ij)\leq \hat{A}(ij)$ for any $i,j\in S$. For example, we set $\hat{A}(ij)\equiv 1$. Therefore, $M(ij)\leq \hat{A}(ij)$ is satisfied. By virtue of Proposition \ref{prop:ext_pote}, there exists $\hat{\ph}\,:\,X_{\hat{A}}\to \R$ such that $[\hat{\ph}]_{k+1}= [\ph]_{k+1}<\infty$, $\hat{\ph}$ is summable, and $\hat{\ph}=\ph$ on $X_{A}$. By Lemma \ref{lem:ex_efunc_geneRuelleop_M_irre_fi}, $\LR_{M,\hat{\ph}}\,:\,F^{k}_{b}(X_{\hat{A}})\to F^{k}_{b}(X_{\hat{A}})$ has the decomposition
$\LR_{M,\hat{\ph}}=\sum_{i=0}^{p-1}\lam_{i}\PR_{i}+\RR$
and $\PR_{i}=h_{i}\otimes \nu_{i}$ for $i$. Proposition \ref{prop:ext_pote} also implies that for any $f\in F^{k}_{b}(X_{A})$, there exists $\hat{f}\in F^{k}_{b}(X_{\hat{A}})$ such that
$f=\hat{f}$ on $X_{A}$ and $\|f\|_{k}=\|\hat{f}\|_{k}$.
Note that such a function $\hat{f}$ is not unique. Remark also that $\LR_{M,\ph} f(\om)=\LR_{M,\hat{\ph}}\hat{f}(\om)$ for $f\in F^{k}_{b}(X_{M})$ if $\om\in X_{A}$. We define operators $\QR_{i}$ and $\SR$ acting on $F^{k}_{b}(X_{A})$ by $\QR_{i} f(\om)=\PR_{i} \hat{f}(\om)$ and $\SR f(\om)=\RR \hat{f}(\om)$ for $f\in F^{k}_{b}(X_{A})$ and $\om\in X_{A}$, respectively. These operators are well-defined. Indeed, if $f_{1},f_{2}$ are in $F^{k}_{b}(X)$ satisfying $f_{1}=f_{2}$ on $X_{A}$, then we see 
$\PR_{i} f_{1}=\nu_{i}(f_{1}\chi_{X_{M}})h_{i}=\nu_{i}(f_{2}\chi_{X_{M}})h_{i}=\PR_{i} f_{2}$ since the support of $\nu_{i}$ is equal to $X_{M}$. Moreover, $\RR f_{1}(\om)=\LR_{M,\hat{\ph}}f_{1}(\om)-\sum_{i=0}^{p-1}\lam_{i}\PR_{i}f_{1}(\om)=\LR_{M,\hat{\ph}}f_{2}(\om)-\sum_{i=0}^{p-1}\lam_{i}\PR_{i}f_{2}(\om)=\RR f_{2}(\om)$ for $\om\in X_{A}$. The equations $\QR_{i}^{2}=\QR_{i}$, $\QR_{i}\QR_{j}=\OR$ for $i\neq j$, $\QR_{i}\LR=\LR\QR_{i}=\lam_{i}\QR_{i}$ and $\QR_{i}\SR=\SR\QR_{i}=\OR$ are valid by the definitions of $\QR_{i}$ and $\RR$. Thus we obtain the spectral decomposition of $\LR_{M,\ph}\,:\,F_{b}^{k}(X_{A})\to F_{b}^{k}(X_{A})$
\ali
{
\textstyle\LR_{M,\ph}=\sum_{i=0}^{p-1}\lam_{i}\QR_{i}+\SR.
}
\indent
Next we show that the spectral radius $\eta$ of $\LR_{M,\ph}$ equals the spectral radius $\lam$ of $\LR_{M,\hat{\ph}}$. 
Letting $g:=h_{0}|_{X_{A}}$, the equation $\LR_{M,\hat{\ph}}h_{0}=\lam h_{0}$ implies $\LR_{M,\ph} g=\lam g$ and then $\lam$ is an eigenvalue of $\LR$. Therefore $\lam\leq \eta$. Moreover when we take $\hat{f}\in F^{k}_{b}(X)$ so that $\|f\|_{k}=\|\hat{f}\|_{k}$, we obtain $\|\LR_{M,\ph}^{n} f\|_{k}\leq \|\LR_{M,\hat{\ph}}^{n}\hat{f}\|_{k}\leq \|\LR_{M,\hat{\ph}}^{n}\|_{k}\|f\|_{k}$. Then $\eta=\lim_{n\to \infty}\|\LR_{M,\ph}^{n}\|_{k}^{1/n}\leq \lim_{n\to \infty}\|\LR_{M,\hat{\ph}}^{n}\|_{k}^{1/n}=\lam$ holds and we get $\eta=\lam$.
\smallskip
\par
We will prove that $\LR_{M,\ph}$ has a spectral gap at $\lam$. By a similar argument above, we see $\|\SR^{n}\|_{k}\leq \|\RR^{n}\|_{k}$ for $n\geq 1$ and thus the spectral radius of $\SR$ is not larger than $r(\RR)<\lam$.
\smallskip
\par
Finally we check the simplicity of $\eta\kappa^{i}=\lam_{i}$. If $f\in F^{k}_{b}(X_{A})$ satisfies $\LR_{M,\ph} f=\lam_{i} f$, then for any $j\neq i$, the equation
$0=\nu_{j}(\LR_{M,\ph} f-\lam_{i}f)=\nu_{j}(\LR_{M,\hat{\ph}}\hat{f}-\lam_{i}\hat{f})=(\lam_{j}-\lam_{i})\nu_{j}(\hat{f})$
holds by $\supp\,\nu_{j}=X_{M}$ and thus $\nu_{j}(\hat{f})=0$. This means $\QR_{j}f=0$. Furthermore, 
$f=\lam_{i}^{-n}\LR^{n} f=\QR_{i}f+\lam_{i}^{-n}\SR^{n}f\to \QR_{i} f$
in $F^{k}_{b}(X_{A})$ as $n\to \infty$ and then $f=\QR_{i}f=\nu_{i}(f) h_{i}$ on $X_{A}$. Hence $\lam_{i}$ is simple. Consequently, we obtain the assertion by replacing $Q_{i}$ and $\SR$ with $\PR_{i}$ and $\RR$, respectively.
\proe
\subsection{Proof of Corollary \ref{cor:ex_efunc_geneRuelleop_M_irre}}\label{sec:proof_cor:ex_efunc_geneRuelleop_M_irre}
Put $S_{2}=S\setminus S_{1}$. Denoted by $\LR_{ij}f=\chi_{S_{i}}\LR(\chi_{S_{j}}f)$ for $i,j=1,2$, where $\chi_{S_{i}}$ is the indicator of $\bigcup_{a\in S_{i}}[a]$. In view of Proposition \ref{prop:ex_efunc_geneRuelleop_noper_ex}, we see $r(\LR_{22})=r(\LR_{M(S_{2})})<\lam$. By Theorem \ref{th:ex_efunc_geneRuelleop_M_irre}, we obtain the spectral decompositions $\LR_{11}=\sum_{i=0}^{p-1}\lam_{i}\PR_{i}+\RR$ and $\PR_{i}=h_{i}\otimes \nu_{i}$. Let $\ti{h}_{i}=(\lam_{i}\IR-\LR_{22})^{-1}\LR_{21}h_{i}$ and $\ti{\nu}_{i}=(\LR_{12}(\lam_{i}\IR-\LR_{22})^{-1})^{*}\nu_{i}$. We consider the decomposition $\LR_{M}
=\sum_{i=0}^{p-1}\lam_{i}\ti{\PR}_{i}+\ti{\RR}$ with
\ali
{
\ti{\PR}_{i}=
\MatII{
h_{i}\otimes \nu_{i}&h_{i}\otimes \ti{\nu}_{i}\\
\ti{h}_{i}\otimes \nu_{i}&\ti{h}_{i}\otimes \ti{\nu}_{i}\\
},\ 
\ti{\RR}=
\MatII
{
\RR&\LR_{12}-\sum_{i=0}^{p-1}\lam_{i}h_{i}\otimes \ti{\nu}_{i}\\
\LR_{21}-\sum_{i=0}^{p-1}\lam_{i}\ti{h}_{i}\otimes \nu_{i}&\LR_{22}-\sum_{i=0}^{p-1}\lam_{i}\ti{h}_{i}\otimes\ti{\nu}_{i}
}
}
by displaying as operator matrix. Note that there is no $M$-admissible word $w\in S^{n}$ such that $w_{1},w_{n}\in S_{i}$ and $w_{l}\in S_{j}$ for some $1<l<n$ for any $\{i,j\}=\{1,2\}$. It follows from this observation in addition to the form $(\lam_{i}\IR-\LR_{22})^{-1}=\sum_{l=0}^{\infty}\LR_{22}^{l}/\lam_{i}^{l+1}$ that $\LR_{12}\LR_{21}=\LR_{21}\LR_{12}=\OR$, $\LR_{21}^{*}\ti{\nu}_{i}=0$, $\LR_{12}\ti{h}_{i}=0$ and $\ti{\nu}_{j}(h_{i})=\nu_{j}(\ti{h}_{i})=\ti{\nu}_{j}(\ti{h}_{i})=0$ for any $i,j$.
Therefore we get 
$(\ti{\PR}_{i})^{2}=\ti{\PR}_{i}$ and $\ti{\PR}_{i}\ti{\RR}=\ti{\RR}\ti{\PR}_{i}=\OR$. Moreover, the simplicity of $\lam_{i}$ follows from $\LR_{M}f=\lam_{i} f$ iff $\LR_{11}(f\chi_{1})=\lam_{i} (f\chi_{1})$ and $f\chi_{2}=(\lam_{i}\IR-\LR_{22})^{-1}\LR_{21}(f\chi_{1})$ iff $f=c(h_{1}+\ti{h}_{1})$ for some $c\in \C$ (see \cite[Proposition 2.7]{T2020} for a proof). Similarity, $\LR_{M}^{*}\mu=\lam_{i}\mu$ iff
$\mu=c(\nu_{i}+\ti{\nu}_{i})$ for some $c\in \C$. By the form $\ti{\PR}_{i}=(h_{i}+\ti{h}_{i})\otimes (\nu_{i}+\ti{\nu}_{i})$, $\ti{\PR}_{i}$ is a projection onto the one-dimensional eigenspace of the eigenvalue $\lam_{i}$.
Finally, we show $r(\ti{\RR})<\lam$. By Proposition \ref{prop:ex_efunc_geneRuelleop_noper_ex}(1), we see $r(\ti{\RR})=\max(r(\RR),r(\RR_{1}))$ with $\RR_{1}=\LR_{22}+\QR$ and $\QR=-\sum_{i=0}^{p-1}\lam_{i}\ti{h}_{i}\otimes \ti{\nu}_{i}$. By the fact $\QR\LR_{22}^{i}\QR=\OR$ for any $i\geq 0$, we obtain the form $\RR_{1}^{n}=\LR_{22}^{n}+\sum_{i=0}^{n-1}\LR_{22}^{i}\QR\LR_{22}^{n-i-1}$ and thus we get $r(\RR_{1})\leq r(\LR_{22})<\lam$.  Hence the proof is complete by replacing $h_{i}:=h_{i}+\ti{h}_{i}$ and $\nu_{i}:=\nu_{i}+\ti{\nu}_{i}$.
\qed
\subsection{Proof of Theorem \ref{th:ex_efunc_geneRuelleop_M}(2)}\label{sec:proof_main2}
\lems\label{lem:ress=0insin_orb}
Assume that (A.1)-(A.3) are satisfied. Assume also that  $X_{A}$ is a finite single orbit, i.e. $X_{A}=\{\om,\si\om,\cdots,\si^{q-1}\om\}$. Then $r_{\text{ess}}(\LR_{M})=0$.
\leme
\pros
Notice that $A$ is finitely irreducible in that case. Moreover, any function $\ph\,:\,X\to \K$ is in $F^{k}_{b}$ for all $\theta\in (0,1)$. Thus  $r_{\text{ess}}(\LR_{M})=0$ from Lemma \ref{lem:ress<=thetalam}.
\proe
\lems\label{lem:ress<lam}
Assume that (A.1)-(A.3) are satisfied. Then $r_{\text{ess}}(\LR_{M})<\lam$.
\leme
\pros
By virtue of Theorem \ref{th:ex_efunc_geneRuelleop_M_irre}, the operator $\LR_{M(T)}$ has spectral gap at $\lam$ and the peripheral eigenvalues of this operator consists of at most a finite number of simple eigenvalues for all $T\in S/\leftrightarrow$. Thus it follows from Proposition \ref{prop:ex_efunc_geneRuelleop_noper_ex}(1)(2) that $\LR_{M}$ has also spectral gap at $\lam$ and the essential spectral radius is less than $\lam$. 
\proe
\lems\label{lem:ress>=lamtheta}
Assume that (A.1)-(A.3) are satisfied. Assume also that $X_{A}$ is not a finite single orbit. Then $r_{\text{ess}}(\LR_{M})\geq \lam\theta$. In particular, for any $r\in [0,\lam\theta)$ except for at most a countable number, $p\in \C$ with $|p|=r$ is an eigenvalue with infinite multiplicity.
\leme
\pros
Then we will prove that $p\in \C$ with $0<|p|<\theta\lam$ is an eigenvalue of $\LR_{M}$ under $X_{A}$ is not finite single orbit. Take the function $g$ of $\LR_{M}$ in Theorem \ref{th:ex_efunc_geneRuelleop_M}. 
We consider the three cases:
\smallskip
\\
Case I: $M$ is irreducible and $X_{M}$ is not finite single orbit. 
In this case, we take a period point $\up\in X_{M}$ with $\up=\si^{l}\up$ and $j\in S_{0}$ so that $M(j\up_{1})=1$, $\up_{0}\neq j$ and $\up_{0},\dots,\up_{l-1}$ are distinct.
Put $f:=(\chi_{[\up_{0}\up_{1}\cdots\up_{k-1}]}-\chi_{[j\up_{1}\cdots\up_{k-1}]})e^{-\ph}$ and $f_{m}:=f(\chi_{[\up_{1}\cdots \up_{ml-1}]}\circ\si)$ for $m\geq k$. Then it is not hard to check that $f_{m}\in F^{k}_{b}(X)$ and $f_{m}\in \ker \LR_{M}$. Recall the eigenfunction $g$ of $\LR_{M}$ given in Theorem \ref{th:ex_efunc_geneRuelleop_M}(1). Put $\ti{g}(\om)=g(\om)$ if $\om\in\supp\,g=\bigcup_{a\in S_{0}}[a]$ and $\ti{g}(\om)=1$ otherwise. Let
$f_{p,m}=g\sum_{n=0}^{\infty}(p/\lam)^{n}(f_{m}/\ti{g})\circ \si^{n}$.
It follows from the fact $f_{m}/\ti{g}\in F^{k}_{b}(X)$ that $f_{p,m}\in F^{k}_{b}(X)$. Notice the equation
\ali
{
\LR_{M}f_{p,m}=&\LR_{M}\Big(g\frac{f_{m}}{\ti{g}}\Big)+\LR_{M}g\Big(\sum_{n=1}^{\infty}\Big(\frac{p}{\lam}\Big)^{n}\Big(\frac{f_{m}}{\ti{g}}\Big)\circ \si^{n-1}\Big)=pf_{p,m}.
}
Now we also check that $\{f_{p,m}\}_{m\geq m_{0}}$ is independent for a large $m_{0}$. By the irreducibly of $M$ and by the assumption $X_{M}\neq \{\up,\si\up,\cdots,\si^{l-1}\up\}$, there exists an $M$-admissible word $w=w_{1}\cdots w_{n}$ such that $M(\up_{s}w_{1})=1$ and $w_{1}\neq \up_{s+1}$ for some $0\leq s<l$ and $M(w_{n}w_{1})=1$. Put $\up^{m}=\up_{0}\cdots \up_{ml-1}\up_{0}\cdots \up_{s}\cd w\cd w\cdots$. Then for any $m$ with $ml> |w|$, we see $f_{p,m}(\up^{m})=f_{m}(\up^{m})=e^{-\ph(\up^{m})}\neq 0$ and $f_{p,m+i}(\up^{m})=0$ for all $i\geq 1$. This says that $\sum_{m\geq |w|/l}c_{m}f_{p,m}=0$ for $c_{m}\in \C$ implies $c_{m}=0$ for all $m$. Thus we obtain the assertion.
\smallskip
\\
Case II: $M$ is irreducible and $X_{M}=\{\up,\si\up,\dots,\si^{l-1}\up\}$ with $\si^{l}\up=\up$. In this case, there exist an $A$-admissible word $w=w_{1}\cdots w_{n}$ and $0\leq s,t<l$ such that $A(\up_{s}w_{1})=A(w_{n}\up_{t})=1$ and $M(\up_{s}w_{1})=M(w_{1}w_{2})=\cdots=M(w_{n-1}w_{n})=0$. 
We put $f_{m}=\chi_{[\up_{s}\cd w\cd \up_{t}\cdots \up_{lm-1}]}$. Then we see $\xi f_{m}\in \ker \LR_{M}$ for any $\xi\in F_{b}^{k}(X)$. Put
$f_{p,m}=g\sum_{n=0}^{\infty}(p/\lam)^{n}f_{m}\circ \si^{n}$.
Then $f_{p,m}\in F^{k}_{b}(X)$ and
$\LR_{M}f_{p,m}=p f_{p,m}$ hold. Now we will show that $\{f_{p,m}\}$ are all nonzero and independent. Put $w^\p=w\cd\up_{t}\cdots \up_{s}$ and 
 $\up^{m}=\up_{s}\cd w\cd\up_{t}\cdots \up_{lm-1}\cd w^\p\cd w^\p\cdots$. Then in addition to the fact $[v_{s}]\subset \supp\, g$, we get $f_{p,m}(\up^{m})=g(\up^{m})f_{m}(\up^{m})\neq 0$ and $f_{p,m+i}(\up^{m})=0$ for all $i\geq 1$. Thus $p$ is an eigenvalue with infinite multiplicity.
\smallskip
\\
Case III: $M$ is not irreducible. We will check that if any $p\in \C$ satisfying that (i) $0<|p|<\theta\lam$, (ii) $|p|\neq r(\LR_{M(T)})\theta$ for any $T\in S/\!\!\leftrightarrow$ and (iii) $p$ is not eigenvalue of $\LR_{M(T)}$ with $|p|>r(\LR_{M(T)})\theta$ for any $T\in S/\!\!\leftrightarrow$, then $p$ is an eigenvalue of $\LR_{M}$. 
 
By the summability of $\ph$ together with Proposition \ref{prop:ex_efunc_geneRuelleop_noper_ex}(1), there exists $S_{1}\in S/\!\!\leftrightarrow$ such that letting $\TR:=\{T\in S/\!\!\leftrightarrow\,:\,S_{1}\preceq T \text{ and }T\neq S_{1}\}$ and $S_{2}:=\bigcup_{T\in \TR}T$, $p$ satisfies $|p|<r(\LR_{M(S_{1})})\theta$ and $|p|>r(\LR_{M(S_{2})})\theta$. Since $p$ is in the resolvent set of $\LR_{M(T)}$ for any $T\in \TR$, it follows from Proposition \ref{prop:ex_efunc_geneRuelleop_noper_ex}(2) that $p$ is in the resolvent of $\LR_{M(S_{2})}$. 
We define $\LR_{ij}f=\chi_{S_{i}}\LR_{M}(\chi_{S_{j}}f)$, where $\chi_{S_{i}}$ denotes the indicator of the set $\Si_{i}=\bigcup_{a\in S_{i}}[a]$.
Let $g_{1}\in F^{k}_{b}(X)$ be an eigenfunction of the eigenvalue $p$ of $\LR_{M(S_{1})}$ and put $g_{2}=(p\IR-\LR_{M(S_{2})})^{-1}\LR_{21}g_{1}$ and $g=g_{1}+g_{2}$. Note that there is no path from a state in $S_{1}\cup S_{2}$ to a state in $S\setminus (S_{1}\cup S_{2})$. By this observation, we obtain
$\LR_{M}g=\LR_{M(S_{1}\cup S_{2})}g=\LR_{M(S_{1})}g_{1}+\LR_{21}g_{1}+\LR_{M(S_{2})}g_{2}=pg$.
Hence we get the assertion. 
\proe
\pros[Proof of Theorem \ref{th:ex_efunc_geneRuelleop_M}(2)] This follows from Lemma \ref{lem:ress=0insin_orb}, Lemma \ref{lem:ress<lam} and Lemma \ref{lem:ress>=lamtheta}.
\proe
\subsection{Proof of Corollary \ref{cor:ex_efunc_geneRuelleop_M_fi}}
It is guaranteed by Lemma \ref{lem:ress<=thetalam} and Lemma \ref{lem:ress>=lamtheta}.
\qed
\subsection{Proof of Theorem \ref{th:escaperate}}\label{sec:escaperate}
First we prove $\limsup_{n\to \infty}\mu_{A}(\Si^{n})^{1/n}\leq \lam_{A}^{-1}r(\LR_{M})$. For $n\geq 1$, $f\in F^{k}_{b}(X)$ and $\om\in X$, we note the equation
\alil
{
\LR_{M}^{n}f(\om)=\sum_{w\in S_{0}^{n}\,:\,w\cd\om_{0}\in W_{n+1}(M)}e^{S_{n}\ph(w\cd\om)}f(w\cd\om)=\LR_{A}^{n}(\chi_{\Si^{n-1}}f)(\om).\label{eq:LM=LASin-1}
}
Therefore
$\mu_{A}(\Si^{n-1})\leq \|h_{A}\|_{\infty}\nu_{A}(\Si^{n-1})=\|h_{A}\|_{\infty}\lam_{A}^{-n}\nu_{A}(\LR_{M}^{n}1)\leq \|h_{A}\|_{\infty}\lam_{A}^{-n}\|\LR_{M}^{n}1\|_{\infty}$.
Thus the assertion holds by using the fact $r(\LR_{M})=r_{C}(\LR_{M})=\lim_{n\to \infty}\|\LR_{M}^{n}1\|_{\infty}^{1/n}$.
\smallskip
\par
Next we show $\liminf_{n\to \infty}\mu_{A}(\Si^{n})^{1/n}\geq \lam_{A}^{-1}r(\LR_{M})$. Recall the quotient space $S/\!\!\leftrightarrow$ defined in Section \ref{sec:Shiftsp}. Choose any $T\in S/\!\!\leftrightarrow$ so that $X_{M(T)}$ has a periodic point. We see
\alil
{
\mu_{A}(\Si^{n-1})=\lam_{A}^{-n}\nu_{A}(\LR_{M}^{n}h_{A})\geq \lam_{A}^{-n}\nu_{A}(\LR_{M(T)}^{n}h_{A})\label{eq:muA>=...}
}
using the equation (\ref{eq:LM=LASin-1}). By virtue of Theorem \ref{th:ex_efunc_geneRuelleop_M_irre} replacing $M$ by $M(T)$, we have the spectral decomposition $\LR_{M(T)}^{n}=\eta^{n}\sum_{i=0}^{p-1}\kappa^{in}\PR_{i}+\RR^{n}$, where we put $\eta=r(\LR_{M(T)})$ and the notation $\kappa, \PR_{i},\RR$ are given in this theorem. For any $n\geq 1$, we denote $q(n)=n \mod p$. We notice the equation
\ali
{
\sum_{i=0}^{p-1}\kappa^{in}\PR_{i}f=\sum_{i,j,l=0}^{p-1}\kappa^{i(n+j-l)}h\chi_{\Si_{j}}\nu(f\chi_{\Si_{l}})=p\sum_{l=0}^{p-1}h\chi_{\Si_{l+q(n)}}\nu(f\chi_{\Si_{l}}).
}
Thus we obtain for any large $n$
\alil
{
\nu_{A}(\LR_{M(T)}^{n}h_{A})=\eta^{n}\nu_{A}(p\sum_{l=0}^{p-1}h\chi_{\Si_{l+q(n)}}\nu(h_{A}\chi_{\Si_{l}}))+\nu_{A}(\RR^{n}h_{A})\geq \frac{c_{\adr{hinui}}}{2}\eta^{n}\label{eq:LRMSi>=...}
} 
with $c_{\adl{hinui}}:=\min_{0\leq q<p}p\sum_{l=0}^{p-1}\nu_{A}(h\chi_{\Si_{l+q}}\nu(h_{A}\chi_{\Si_{l}}))>0$, where the last inequality uses the fact $\|\eta^{-n}\RR^{n}h_{A}\|_{\infty}\leq c_{\adr{hinui}}/2$ for any large $n$. It follows from (\ref{eq:muA>=...}) and (\ref{eq:LRMSi>=...}) that $\liminf_{n\to \infty}\mu_{A}(\Si^{n})^{1/n}\geq \lam_{A}^{-1}r(\LR_{M(T)})$ for any $T\in S/\!\!\leftrightarrow$. Proposition \ref{prop:ex_efunc_geneRuelleop_noper_ex}(1) says $r(\LR_{M})=\max_{T\in S/\leftrightarrow}r(\LR_{M(T)})$ and consequently we get $\liminf_{n\to \infty}\mu_{A}(\Si^{n})^{1/n}\geq \lam_{A}^{-1}r(\LR_{M})$. Hence the assertion holds together with the facts $\lam=\exp(P(\ph|_{X_{M}}))$ an $\lam_{A}=\exp(P(\ph))$.
\qed
\section{Applications and examples}\label{sec:app_ex}
\subsection{Convergence of the topological pressure and the Gibbs measure of perturbed potential in our open system setting}\label{sec:conv_pres_Gibbs}
\prop
{\label{prop:conv_pres}
Assume that the conditions (A.1)-(A.3) and (B.1)-(B.3) are satisfied. Then the topological pressure of $\phe$ converges to the topological pressure of $\ph|_{X_{M}}$.
}
\pros
Choose any limit point $\eta$ of $\{\exp(P(\phe))\}_{\e}$.
By Proposition \ref{prop:conv_efunc_LB}, $\eta$ becomes an eigenvalue of $\LR_{M}$ and therefore $\eta\leq \lam:=r(\LR_{M})$. On the other hand, we can choose $T\in S/\!\!\leftrightarrow$ so that $\exp(P(\ph|_{X_{M(T)}}))=\lam$ together with 
Proposition \ref{prop:ex_efunc_geneRuelleop_noper_ex}(1). Let $N=\bigcup_{M(ij)=1, M(T)(ij)=0}[ij]$ and $\psi(\e,\cd)=\ph(\e,\cd)-(1/\e)\chi_{N}$.
Let $\lam_{0}$ be any limit point of $\{\exp(P(\psi(\e,\cd)))\}_{\e}$. Then we may assume $\lam_{0}\leq \eta$ by the definition of $\psi(\e,\cd)$. We will show $\lam_{0}=\lam$.
Take a nonnegative eigenfunction $g_{0}$ of $\lam_{0}$ of $\LR_{M(T)}$ from Proposition \ref{prop:conv_efunc_LB}, and a positive eigenvector $\nu$ of $\lam$ of $\LR_{M(T)}^{*}$ from Theorem \ref{th:ex_efunc_geneRuelleop_M_irre}. Since $M(T)$ is irreducible, $\nu(g_{0})>0$ and therefore the equation $(\lam_{0}-\lam)\nu(g_{0})=\nu((\LR_{M(T)}-\LR_{M(T)})g_{0})=0$ implies $\lam_{0}=\lam$. Thus $\lam=\eta$. Hence the proof is complete.
\proe
\prop
{\label{prop:conv_Gibbs}
Assume that the conditions (A.1)-(A.3) and (B.1)-(B.3) are satisfied. Assume also that $\{T\in S/\!\!\leftrightarrow\,:\,P(\ph|_{X_{M(T)}})=P(\ph|_{X_{M}})\}$ is only one element $T$. Take the pair $(h,\nu)$ in Corollary \ref{cor:ex_efunc_geneRuelleop_M_irre}.
Then the Gibbs measure $\mu(\e,\cd)$ of $\phe$ converges to the $\si$-invariant Borel probability measure $h\nu$ weakly as $\e\to 0$.
}
\pros
Take $(\lam(\e),h(\e,\cd),\nu(\e,\cd))$ as the spectral triplet given as (\ref{eq:pTSC}) and $(\lam,h_{0},\nu_{0})$ as the spectral triplet given in Corollary \ref{cor:ex_efunc_geneRuelleop_M_irre}. Put $g(\e,\cd)=h(\e,\cd)/\|h(\e,\cd)\|_{\infty}$. Note the form $h(\e,\cd)=g(\e,\cd)/\nu(\e,g(\e,\cd))$. First we state $\nu(\e,\cd)\to \nu_{0}$ weakly. This is yielded by the existence of limit point of $\nu(\e,\cd)$ (Proposition \ref{prop:conv_evec_LB}) and by the simplicity of $\lam$ of the dual $\LR_{M}^{*}$ (Corollary \ref{cor:ex_efunc_geneRuelleop_M_irre}).
Next we show $\sup_{\e>0}\|h(\e,\cd)\|_{k}<\infty$. To do this, we check $\liminf _{\e\to 0}\nu(\e,g(\e,\cd))>0$. Fix an $M$-admissible word $w\in T^{k}$. Since $g(\e,\up)\leq e^{c_{\adr{udphe2}}}g(\e,\om)$ holds for any $\om,\up\in [w]$ and since $g(\e,\om)$ has a limit point $c_{\adr{pg}}h_{0}(\om)>0$ for some $c_{\adl{pg}}>0$, we obtain $\liminf_{\e\to 0}\inf_{[w]}g(\e,\cd)\geq e^{-c_{\adr{udphe2}}}c_{\adr{pg}}h_{0}(\om)>0$. Therefore $\nu(\e,g(\e,\cd))\geq \nu(\e,[w])\inf_{[w]}g(\e,\cd)$ and the right hand side has a positive lower bound by $\nu(\e,[w])\to \nu_{0}([w])>0$. Consequently $\{h(\e,\cd)\}$ has a finite upper bound and satisfies $\limsup_{\e\to 0}\|h(\e,\cd)\|_{k}<\infty$.
Finally we prove convergence of $\mu(\e,\cd)=h(\e,\cd)\nu(\e,\cd)$. Let $f\in F_{b}^{k}(X)$. By using Theorem \ref{th:conv_evec}, we have the estimate of $\kappa(\e,f):=\nu(\e,h(\e,\cd)f)/\nu(\e,h_{0})$:
\ali
{
|\kappa(\e,f)-\nu(h(\e,\cd)f)|=&|\kappa(\e,(\LR_{A,\phe}-\LR_{M})(h\otimes \kappa(\e,\cd)-\IR)(\ER-\lam\IR)^{-1}(h(\e,\cd)f))|\\
\leq& \frac{\|\LR_{A,\phe}-\LR_{M}\|_{\infty}}{\nu(\e,h_{0})}\Big(\frac{\|h\|_{\infty}}{\nu(\e,h_{0})}+1\Big)\|(\ER-\lam\IR)^{-1}\|_{k}\|h(\e,\cd)f\|_{k}\to 0
}
as $\e\to 0$ by using convergence $\nu(\e,h_{0})\to \nu_{0}(h_{0})=1$ and Proposition \ref{prop:conv_Ruelleop_LB}, where $\ER:=\LR_{M}-\lam(h_{0}\otimes \nu_{0})$. Moreover, when we take a limit point $c_{\adl{ph}}h_{0}$ of $\{h(\e,\cd)\}$ with a constant $c_{\adr{ph}}>0$, we see $\nu(h(\e,\cd)f)\to c_{\adr{ph}}\nu_{0}(h_{0}f)$ as $\e\to 0$ running through a suitable sequence. If we put $f=1$, then we see that $c_{\adr{ph}}$ must be $1$ by $\kappa(\e,1)\to 1$. Thus $\nu(h(\e,\cd)f)\to \nu_{0}(h_{0}f)$. Hence the assertion is valid by noting $\nu_{0}(h_{0} f)=\nu_{0}(\chi_{\bigcup_{a\in T}[a]}h_{0} f)=\nu(h f)$.
\proe
\subsection{Example (A Markov measure with countable states)}\label{sec:ex}
Put $S=\{1,2,\dots\}$ and
take positive numbers $a_{n}, b_{n}$ with $\sum_{n\in S}\max\{a_{n},b_{n}\}<+\infty$. Let $S\times S$ zero-one matrix $M=(M(ij))$ be $M(ij)=1$ if $j=1$ or $j=i+1$, and $M(ij)=0$ otherwise.
We set $\ph\,:\,X_{M}\to \R$ by $\ph(\om)=\log a_{\om_{0}}$ if $\om_{1}=1$ and $\ph(\om)=\log b_{\om_{0}}$ if $\om_{1}=\om_{0}+1$. Then $\ph$ is summable.
The below follows immediately from Theorem \ref{th:ex_efunc_geneRuelleop_M_irre} with $k=1$.
\props
\label{prop:ex_prop1}
Under the above notation, the Ruelle operator $\LR\,:\,F^{1}_{b}(X_{M})\to F^{1}_{b}(X_{M})$ of $\ph$ has the decomposition $\LR=\lam (h\otimes \nu)+\RR$ and the spectral radius of $\RR$ is less than $\lam$, where $(\lam,h,\nu)$ appears in Theorem \ref{th:ex_efunc_geneRuelleop_M_irre}. 
\prope
We define a stochastic matrix $P=(P(ij))$ indexed by $S\times S$ as $P(i1)=b_{1}b_{2}\cdots b_{i-1}a_{i}/\lam^{i}$, $P(i\ i+1)=1$ and $P(ij)=0$ otherwise. Then the measure $\mu=h\nu$ becomes the Markov measure of $P$ which runs backwards. We can check that the potential $\om\mapsto \log P(\om_{0}\om_{1})$ is cohomologous to $\ph-P(\ph)$ via the transfer function $\log h$. 
\subsection{Example (Graph iterated function systems)}\label{sec:ex2}
Let $G=(V,E)$ be a directed multigraph with the countable set $V$ of vertices and the countable set $E$ of edges. For $e\in E$, denoted by $i(e)$ the initial vertex of $e$ and by $t(e)$ the terminal vertex of $e$. We introduce a set $(G,(J_{v})_{v\in V},(O_{v})_{v\in V}, (T_{e})_{e\in E})$ as follows:
\ite
{
\item[(G.1)] The graph $G$ is strongly connected, i.e. there is a path between any two vertices.
\item[(G.2)] Each $J_{v}$ is a nonempty closed subset of a Banach space $(Y_{v},\|\cd\|_{v})$ with $c_{\adl{dJv}}:=\sup_{v}\diam J_{v}<+\infty$.
\item[(G.3)] Each $O_{v}$ is an open convex subset of $Y_{v}$ and $J_{v}\subset O_{v}$.
\item[(G.4)] Each $T_{e}$ is a $C^{1}$ map from $O_{t(e)}$ to $O_{i(e)}$ satisfying  $T_{e}J_{t(e)}\subset J_{i(e)}$. There exists $0<r<1$ such that the norm $\|DT_{e}(x)\|$ of the derivative of $T_{e}$ at $x$ is no larger than $r$ for any $e\in E$ and $x\in O_{t(e)}$. Moreover, there exist constants $c_{\adl{DTe}}>0$ and $0<\beta\leq 1$ such that $|\|DT_{e}(x)\|-\|DT_{e}(y)\||\leq c_{\adr{DTe}}\|DT_{e}(x)\|\|x-y\|_{t(e)}^{\beta}$ for any $x,y\in J_{t(e)}$ and $e\in E$.
\item[(G.5)] $s_{*}:=\inf\{s\geq 0\,:\,\sum_{e\in E}\sup_{x\in J_{t(e)}}\|DT_{e}(x)\|^{s}<\infty\}$ is defined as an real number.
}
Let $M=(M(ee^\p))$ be $E\times E$ zero-one matrix defined by $M(ee^\p)=1$ if $t(e)=i(e^\p)$ and $M(ee^\p)=0$ otherwise. 
Then $X=X_{M}$ is topological transitive by (G.1). For $\om\in X$, we define $\pi\om\in Y_{i(\om_{0})}$ by
$\pi\om=\bigcap_{n=0}^{\infty}T_{\om_{0}}\circ\cdots\circ T_{\om_{n}}(J_{t(\om_{n})})$.
This is well-defined by $\|DT_{e}\|\leq r<1$ 
and is called the coding map of the limit set $\pi(X)$. Let $\ph\,:\,X\to \R$ be
\alil
{
\ph(\om)=\log \|DT_{\om_{0}}(\pi\si\om)\|.\label{eq:phys_po}
}
Then the thermodynamic future of the potential $s\ph$ ($s>s_{*}$) provides us important information in the fractal analysis. In fact, under a suitable condition for $(G,(J_{v}),(O_{b}), (T_{e}))$, the Hausdorff dimension $\dim_{H}\pi(X)$ of the limit set $\pi(X)$ is given by $\inf\{s>0\,:\,P(s\ph)<0\}$ which is so-called a generalized Bowen formula (e.g. \cite{MSU,MU,PU, NPL, Priyadarshi}). Moreover, the spectral decomposition of a (unperturbed) Ruelle operator $\LR_{s\ph}$ of $s\ph$ is an important role for perturbation analysis of $(G,(J_{v}),(O_{b}), (T_{e}))$ \cite{T2020pre, T2019, T2016}. In previous work, it is known in \cite{MU2001,MU, Priyadarshi} that if the incidence matrix $M$ is finitely primitive, then $\LR_{s\ph}$ has a gap property in the eigenvalue $\exp(P(s\ph))$. By using results in present paper, it turns out that the finite irreducibility is not needed to obtain such a gap property:
\thm
{\label{th:GIFS}
Assume that the conditions (G.1)-(G.5) are satisfied. Take the function $\ph$ of (\ref{eq:phys_po}) and $\theta=r^{\beta}$. Then the Ruelle operator of the potential $s\ph$ acting on $F_{b}^{1}(X_{M})$ has the decomposition (\ref{eq:LMph=sumi=0p-1...}) replacing by $\ph:=s\ph$ for each $s>s_{*}$.
}
\pros
It is sufficient to check the conditions of Theorem \ref{th:ex_efunc_geneRuelleop_M_irre}.
Since $G$ is strongly connected, the incidence matrix $M$ is irreducible. The summability of $s\ph$ follows from the condition (G.5) and $\|DT_{e}\|\leq 1$. Further, the condition (G.4) together with the chain rule for the derivatives of $T_{e}$ implies that for any $\om_{0}\cdots\om_{n-1}=\up_{0}\cdots\up_{n-1}$ and $\om_{n}\neq \up_{n}$ with $n\geq 1$, $|s\ph(\om)-s\ph(\up)|\leq s c_{\adr{dJv}}c_{\adr{DTe}}r^{(n-1)\beta}$. Thus $\ph$ is a locally $r^\beta$-Lipschitz continuous function. Hence we obtain the assertion by applying Theorem \ref{th:ex_efunc_geneRuelleop_M_irre} to the potential $s\ph$.
\proe
\subsection{Example (locally constant potentials)}\label{sec:ex3}
A function $\ph\,:\,X\to \R$ is called {\it uniformly locally constant} if $[\ph]_{n}=0$ for some $n$ (\cite{Kempton} for terminology). For example, the countable Markov chains \cite{GS, Salama} are of the case $[\ph]_{2}=0$. 
Note that this notion does not depend on choosing $\theta$.
The following immediately follows from Theorem \ref{th:ex_efunc_geneRuelleop_M}:
\thm
{
Assume that the conditions $(A.1)$-$(A.3)$ with fixed integer $k\geq 1$ are satisfied and the potential $\ph$ is uniformly locally constant with $[\ph]_{k+1}=0$. Then the eigenfunction $g$ in Theorem \ref{th:ex_efunc_geneRuelleop_M}(1) of the eigenvalue $\exp(P(\ph))$ of the operator $\LR_{M}\,:\,F_{b}^{k}(X)\to F_{b}^{k}(X)$ is uniformly locally constant with $[g]_{k}=0$. Moreover, for any $\theta\in (0,1)$, the essential spectral radius of $\LR_{M}$ acting on $F_{b}^{k}(X)$ with the metric $d_{\theta}$ is equal to $\theta \exp(P(\ph))$.
}
\appendix
\section{Extension of potentials}\label{sec:Ext_poten}
Let $X=X_{A}$ be a nonempty topological Markov shift with transition matrix $A$ and with countable states $S$. Take a subsystem $Y=X_{M}$ of $X$ with $S\times S$ transition matrix $M$ with $M(ij)\leq A(ij)$ and with $[a]\cap Y\neq \emptyset$ for any $a\in S$. Let $\ph\,:\,Y\to \R$ be a summable function. For each $A$-admissible word $w\in \bigcup_{n=1}^{\infty}S^{n}$ with $[w]\cap Y\neq \emptyset$, fix an element $\nu^{w}\in [w]\cap Y$. We define a function $\hat{\ph}$ on $X$ by
\ali
{
\hat{\ph}(\om)=&
\case
{
\ph(\up^{w}),&\text{ if }w=\om_{0}\cdots \om_{m}\text{ satisfies }[w]\cap Y\neq \emptyset, [w\om_{m+1}]\cap Y=\emptyset \text{ for an }m\geq 0\\
\ph(\om),&\text{ if }\om\in Y.
}
}
\prop
{\label{prop:ext_pote}
Under the above conditions, 
(i) $\hat{\ph}=\ph$ on $Y$; (ii) $\sup_{\om\in [s]}\hat{\ph}(\om)\leq \sup_{\om\in [s]}\ph(\om)$, in particular, $\hat{\ph}$ is summable; (iii) $\|\hat{\ph}\|_{\infty}=\|\ph\|_{\infty}$; (iv) $[\hat{\ph}]_{n}=[\ph]_{n}$ for all $n\geq 1$.
}
\pros
By the definition, $\hat{\ph}$ is an extension of $\ph$ and satisfies $\|\hat{\ph}\|_{\infty}=\|\ph\|_{\infty}$. Moreover, by $\sum_{s\in S}e^{\sup_{\om\in [s]}\hat{\ph}(\om)}\leq \sum_{s\in S}e^{\sup_{\om\in [s]}\ph(\om)}<\infty$, 
$\hat{\ph}$ is also summable.
Now we check (iv). Let $\om,\up \in X$ with $d_{\theta}(\om,\up)=\theta^{n+1}$ for some $n\geq 0$. 
Let $m_{0}\geq 0$ be the smallest number so that $w_{0}:=\om_{0}\cdots \om_{m_{0}}$ satisfies $[w_{0}]\cap Y\neq \emptyset$ and $[w_{0}\om_{m+1}]\cap Y= \emptyset$ if $\om\notin Y$ and $m_{0}=+\infty$ if $\om\in Y$. Similarity, we take $m_{1}\geq 0$ and $w_{1}=\up_{0}\cdots \up_{m_{1}}$ for $\up$.
Consider the four cases:
\smallskip
\\
Case I: $m_{0}<n$. In this case, we see $n_{0}=m_{0}$ and $w_{0}=w_{1}$, and therefore
$|\hat{\ph}(\om)-\hat{\ph}(\up)|=|\ph(\up^{w_{0}})-\ph(\up^{w_{0}})|=0<d_{\theta}(\om,\up)$.
\smallskip
\\
Case II:  $n\leq m_{0}<\infty$ and $m_{1}<\infty$. We notice $\up^{w_{0}},\ \up^{w_{1}}\in [\om_{0}\cdots \om_{n}]$ and thus
$|\hat{\ph}(\om)-\hat{\ph}(\up)|=|\ph(\up^{w_{0}})-\ph(\up^{w_{1}})|
\leq [\ph]_{n+1}\theta^{n+1}=[\ph]_{n+1}d_{\theta}(\om,\up)$.
\smallskip
\\
Case III:  $n\leq m_{0}<\infty$ and $m_{1}=\infty$. By $\up^{w_{0}},\ \up\in [\om_{0}\cdots \om_{n}]$, we get
$|\hat{\ph}(\om)-\hat{\ph}(\up)|=|\ph(\up^{w_{0}})-\ph(\up)|
\leq [\ph]_{n+1}\theta^{n+1}=[\ph]_{n+1}d_{\theta}(\om,\up)$.
\smallskip
\\
Case IV: $m_{0}=\infty$. By switching $\om$ and $\up$, similar arguments above Case I and Case III imply $[\hat{\ph}]_{n+1}\leq [\ph]_{n+1}$ for all $n\geq 0$. 

On the other hand, for $\om,\up\in Y$ with $\om\in [\up_{0}\cdots\up_{n-1}]$, we have
$|\ph(\om)-\ph(\up)|=|\hat{\ph}(\om)-\hat{\ph}(\up)|\leq [\hat{\ph}]_{n}d_{\theta}(\om,\up)$.
Thus $[\ph]_{n}\leq [\hat{\ph}]_{n}$. Hence the proof is complete.
\proe
\section{Perturbation of eigenvectors of linear operators}\label{sec:per_evec}
In this section, we recall a special case given in \cite[Theorem 2.2]{T2022_pre}.
Assume the following:
\ite
{
\item[(I)] $\LR$ is a linear operator acting on a linear space $\BR$ over $\K$.
\item[(II)] There exists a triplet $(\lam,h,\nu)\in \K\times \BR\times \BR^{*}$ such that $\lam\neq 0$, $\LR h=\lam h$, $\LR^{*}\nu=\lam \nu$ and $\nu(h)=1$, where $\BR^{*}$ is the dual space of $\BR$ and $\LR^{*}$ is the dual operator of $\LR$. Moreover, $\ER:=\LR-\lam(h\otimes \nu)$ satisfies that $(\lam\IR-\ER)^{-1}\,:\,\BR\to \BR$ is well-defined as linear operator.
\item[(III)] $\LR(\e,\cd)$ is a parametrized linear operator acting on $\BR$ with $\e>0$ and there exists a pair $(\lam(\e),\nu(\e,\cd))\in \K\times \BR^{*}$ such that $\LR(\e,\cd)^{*}\nu(\e,\cd)=\lam(\e)\nu(\e,\cd)$ and $\nu(\e,h)\neq 0$ for any $\e$.
}
\thm
{\label{th:conv_evec}
Under the conditions above, $\kappa(\e,\cd):=\nu(\e,\cd)/\nu(\e,h)$ has the form
$\kappa(\e,f)=\nu(f)+\kappa(\e,\tLR(\e,\cd)(h\otimes \kappa(\e,\cd)-\IR)(\ER-\lam\IR)^{-1}f)$
for each $f\in \BR$, where $\tLR(\e,\cd):=\LR(\e,\cd)-\LR$.
}
\pros
We start with the equation $\IR-h\otimes \nu=(\LR-\lam)(\ER-\lam\IR)^{-1}$. This is obtained by direct checking $(\IR-h\otimes \nu)(\ER-\lam\IR)=(\ER-\lam\IR)(\IR-h\otimes \nu)=\LR-\lam\IR$.
Note also the equation $\nu(\e,(\lam(\e)-\lam)h)=\nu(\e,(\LR(\e,\cd)-\LR)h)$ using the conditions (II) and (III). 
Moreover, 
\ali
{
\kappa(\e,f)-\nu(f)=\kappa(\e,(\IR-(h\otimes \nu))f)=&\kappa(\e,(\LR-\lam\IR)(\ER-\lam\IR)^{-1}f)\\
=&\kappa(\e,(-\tLR(\e,\cd)+(\lam(\e)-\lam)\IR)(\ER-\lam\IR)^{-1}f).
}
Hence the assertion holds together with $\lam(\e)-\lam=\kappa(\e,\tLR(\e,\cd)h)$.
\proe


\begin{thebibliography}{99}
\bibitem{A} Aaronson, J.: An Introduction to Infinite Ergodic Theory. Am. Math. Soc., Providence (1997).
\bibitem{AD} Aaronson, J., Denker, M.: Local limit theorems for partial sums of stationary sequences generated by Gibbs-Markov maps. Stoch. Dyn. \textbf{1}, no. 2, 193237 (2001).
\bibitem{ADU} Aaronson, J., Denker, M., Urba\'nski, M.: Ergodic theory for Markov
fibred systems and parabolic rational maps, Trans. Amer. Math. Soc. 337, no. 2,
495--548, (1993).
\bibitem{AJS} Ammar, A., Jeribi, A., Saadaoui, B.; Frobenius-Schur factorization for multivalued $2 \times 2$ matrices linear operator, Mediterr. J. Math. \textbf{14} (1), Paper No. 29, 29 pp, (2017).
\bibitem{Baladi} Baladi, V.: Positive transfer operators and decay of correlations,
World Scientific, (2000).
\bibitem{Cyr_Sarig} Cyr, V., Sarig, O.: Spectral gap and transience for Ruelle operators on countable Markov shifts. Comm. Math. Phys., 292, 637--666, (2009).
\bibitem{Demers_} Demers, M., Ianzano, C., Mayer, P., Morfe, P., Yoo, E.: Limiting distributions for countable state topological Markov chains with holes, Discrete and Contin. Dynam. Sys. 37, 105--130, (2017).
\bibitem{GS} Gurewic, B. M., Savchenko, S.V.: Thermodynamic formalism for countable symbolic Markov chains, Russian Math. Surv. 53, 245--344 (1998).
\bibitem{HH} Hennion, H., Hevr\'e, L.: Limit theorems for Markov chains and stochastic properties of dynamical systems by quasi-compactness, 1766, Lectures Notes in Mathematics, Springer-Verlag, Berlin, (2001).
\bibitem{Kempton} Kempton, T.: Zero temperature limits of Gibbs equilibrium states for countable Markov shifts. J. Stat. Phys. 143(4), 795--806 (2011). 
\bibitem{Lin} Lin, M.: Mixing for Markov operators, Z. Wahrscheinlikeit u.v. Gebiete 19, 231--243, (1971).
\bibitem{MSU} Mauldin, R. D., Szarek, T., Urba\'nski, M.: Graph directed Markov systems on Hilbert spaces, Math. Proc. Cambridge Philos. Soc. 147(2), 455--488, (2009).
\bibitem{MU2001} Mauldin, R. D.,  Urba\'nski, M.: Gibbs states on the symbolic space over an infinite alphabet. Isr. J. Math. 125, 93-130 (2001).
\bibitem{MU} Mauldin, R. D., Urba\'nski, M.: Graph Directed Markov Systems : Geometry and dynamics of limit sets, Cambridge (2003).
\bibitem{NPL} R. D. Nussbaum, A. Priyadarshi, and S. Verduyn Lunel, Positive operators and Hausdorff dimension of invariant sets. Trans. Amer. Math. Soc. \textbf{364}, no. 2, 1029--1066 (2012).
\bibitem{PU} Pollicott, M., Urba\'nski, M.: Open Conformal Systems and Perturbations of Transfer Operators, Lecture Notes in Mathematics, Vol. 2206,  Springer, Berlin, (2018).
\bibitem{Priyadarshi} Priyadarshi, A. Infinite Graph-Directed Systems and Hausdorff Dimension. Waves Wavelets Fractals Adv. Anal. 3:84--95 (2017).
\bibitem{Sar09} Sarig, Omri M.: Lecture Notes on Thermodynamic Formalism for Topological Markov Shifts, Penn State (2009).
\bibitem{Sar03} Sarig, Omri M.: Existence of Gibbs measures for countable Markov shifts, Proc. Amer. Math. Soc. \textbf{131}, no. 6, 1751-1758 (2003).
\bibitem{Sar01_2} Sarig, Omri M.: Thermodynamic formalism for null recurrent potentials. Israel J. Math. \textbf{121}, 285--311 (2001).
\bibitem{RSU} Roy, M., Sumi, H., Urba\'nski, M.: Analytic families of holomorphic iterated function systems, Nonlinearity 21, no. 10, 2255--2279, (2008).
\bibitem {Salama} Salama, I. A., Topological entropy and recurrence of countable chains. Pac. J. Math. 134(2), 325--341 (1988).
\bibitem{SU} Stratmann, B. O., Urba\'nski, M.: Pseudo-Markov systems and infinitely generated Schottky groups. Amer. J. Math. 129, no. 4, 1019--1062, (2007).
\bibitem{T2022_pre} Tanaka, H.: General asymptotic perturbation theory in transfer operators, submitted (arXiv:2205.12561).
\bibitem{T2020pre} H. Tanaka, Asymptotic solution of Bowen equation for perturbed potentials on shift spaces with countable states, to appear in J. Fractal Geom.
\bibitem{T2020} Tanaka, H.: Perturbed finite-state Markov systems with holes and Perron complements of Ruelle operators, Israel J. of Math. 236, 91--131, (2020).
\bibitem{T2019} Tanaka, H.: Perturbation analysis in thermodynamics using matrix representations of Ruelle operators and its application to graph IFS, Nonlinearity, 32, 728--767, (2019). 
\bibitem{T2016} Tanaka, H.: Asymptotic perturbation of graph iterated function systems, Journal of Fractal Geometry, {\bf 3} (2016), 119--161.
\bibitem{T2009} Tanaka, H.: Spectral properties of a class of generalized Ruelle operators. Hiroshima Math. J. {\bf 39}, 181--205 (2009).
\bibitem{Tretter} Tretter, C.: Spectral Theory of Block Operator Matrices and Applications. Imperial College Press, London (2008).
\end{thebibliography}
\end{document}